\documentclass[11pt]{article}

\usepackage{epsfig}
\usepackage{epstopdf}
\usepackage{latexsym}
\usepackage{float}
\usepackage{graphicx}
\usepackage{amsfonts, amsbsy}
\usepackage{amsmath, subfigure}
\usepackage{amscd}
\usepackage{hyperref}
\usepackage[colorinlistoftodos]{todonotes}
\usepackage{authblk}

\usepackage{gensymb}

\title{ A bridge between  invariant dynamical structures and uncertainty  quantification}

 \author[1]{G. Garc\'ia-S\'anchez}
 \author[1]{A. M. Mancho}
 \author[2]{S. Wiggins}
 
 \affil[1]{Instituto de Ciencias Matem\'aticas, CSIC, C/Nicol\'as Cabrera 15, Campus Cantoblanco, 28049 Madrid, Spain}
\affil[2]{School of Mathematics, Fry Building, 
Woodland Road, University of Bristol, 
 Bristol BS8 1UG, United Kingdom}

\begin{document}

\maketitle

\begin{abstract}
We develop a new quantifier for forward time uncertainty for trajectories that are solutions of models generated from data sets. Our uncertainty quantifier is defined on the phase space in which the trajectories evolve and we show that it has a rich structure that is directly related to phase space structures from dynamical systems theory, such as hyperbolic trajectories and their stable and unstable manifolds. We apply our approach to an ocean data set, as well as standard benchmark models from deterministic dynamical systems theory. A significant application of our results, is that they allow  a quantitative comparison of the transport performance  described from different ocean data sets. This is particularly interesting nowadays when a wide variety of sources are available, since our methodology provides avenues for assessing the effective use of these data sets in a variety of situations.  
\end{abstract}

\section{Introduction}
Uncertainty Quantification (UQ) searches for a quantitative characterization of uncertainties in computational models of real-world applications. Uncertainties associated to these models have been studied for many years. Depending on the context they can arise as a result of noisy or incomplete data,  errors introduced by numerical models, or by the lack of complete understanding of the governing physical processes. 

In geophysical contexts, uncertainty is a topic of much interest because 
it is inherent to the equations describing the motion of fluids such as the ocean or the atmosphere.  The motivation of this  paper is  trying to acquire a deeper understanding of uncertainty quantification, that allows a better characterization of its presence in ocean models. Large amount of oceanic data are becoming now available. For instance, the Copernicus Marine Environment Monitoring Service (CMEMS) provides regular and systematic information about the physical state and dynamics of the ocean for the global ocean and the European regional seas. The data cover the current and future state of variables and the provision of retrospective data records (re-analysis). Other oceanic services include ocean currents supplied by altimeter satellites, like AVISO; the HYbrid Coordinate Ocean Model (HYCOM), a consortium, which is a multi-institutional effort sponsored by the National Ocean Partnership Program (NOPP), as part of the U. S. Global Ocean Data Assimilation Experiment (GODAE), to develop and evaluate a data-assimilative hybrid coordinate ocean models, etc. 

Ocean models are built on partial differential equations, which are only an approximation of reality. Indeed, typically ocean models approach ocean motion by the Reynolds-averaged Navier-Stokes equations using the hydrostatic and Boussinesq assumptions. In this context inaccurate or inadequate models lead to {\em structural uncertainty}. Even if these equations were perfectly accurate models, they contain parameters, such as viscosity or diffusivity,  which are not known precisely and cause {\em parameter uncertainty}.  Boundary conditions on input variables, lead to {\em parametric variability}. For instance input variables 
such as velocities, experience at the surface atmospheric winds, which impose forcings on the velocities, and these atmospheric forcings are not precisely known. Eventually these models are not solved exactly, but numerically through discretization methods and numerical approximations, which are subjected to numerical errors, and causes {\em algorithmic uncertainty}.
For all these reasons, the predicted ocean currents are also uncertain. Uncertainties in the solutions, system outputs, due to the uncertainties in the system inputs  is referred to as {\em forward uncertainty quantification} \cite{sullivan2015}. 

\begin{figure}[htb!]
  \begin{center}
 \includegraphics[scale=0.7]{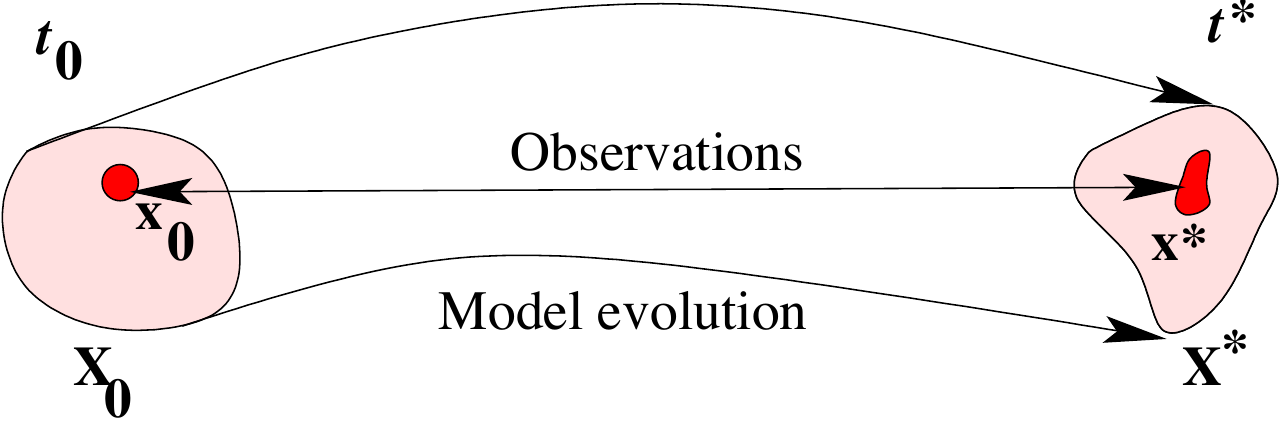} 
  \end{center}
  \caption{A graphical representation of two sequential observations and their evolution according to a mathematical model. The initial observation at time $t_0$ is expressed by the red initial condition ${\bf x_0}$. The final observed state  
 ${\bf x^*}$ at time $t^*$ is referred to as the  "target". The  evolution law for these observations is unknown, but is approached by a  model, which in our setting involves the velocities ${\bf v}({\bf x},t)$ of the system \eqref{v(t)}. The evolution according to a model, of  a neighbourhood of points, ${\bf X_0}$, close to the initial observation, ${\bf x_0}$, is illustrated with the  pink color. }
  \label{fig:uqsk}
\end{figure}

In this work, we are interested in evaluating the reliability of the outputs, ocean currents, regarding the transport they produce. This is a perspective slightly different to the classical one now described, in which ocean models adequacy is judged just against the velocity fields. 
Transport in the  ocean surface is produced by fluids  parcels that follow trajectories ${\bf x}(t)$ that evolve according to the dynamical system:
\begin{equation}
    \frac{d{\bf x}}{dt} = {\bf v}({\bf x},t),
    \label{v(t)}
\end{equation}
where the position is described in longitude ($\lambda$) and latitude ($\phi$) coordinates, that is, ${\bf x}=(\lambda,\phi)$, and ${\bf v}$ represents the velocity field, which has two components determined by the zonal ($u$) and meridional ($v$) velocities. In longitude/latitude coordinates, the dynamical system in Eq.\eqref{v(t)} can be rewritten as:
\begin{equation}
    \begin{cases}
    \dfrac{d \lambda}{d t} = \dfrac{u(\lambda,\phi,t)}{R\cos \phi} \\[.3cm]
    \dfrac{d \phi}{d t} = \dfrac{v(\lambda,\phi,t)}{R}
    \end{cases}\label{v(t)s}
\end{equation}
where $R$ is the Earth's radius. The system \eqref{v(t)}, which is a general expression  encompassing the specific problem dealt with in \eqref{v(t)s}, is a nonlinear non-autonomous dynamical system in 2D. Uncertainties in the velocities ${\bf v}({\bf x},t)$, or more generally in the vector field, produce uncertainties in the solutions ${\bf x}(t)$.
In general, as illustrated in figure  \ref{fig:uqsk},  the exact model connecting two successive observations is not known. Only the initial observation, ${\bf x_0}$,  and the final state,  ${\bf x^*}$, are measurable. They are presented in red color in the figure. The evolution of initial conditions, ${\bf X_0}$, in a neighbourhood  close to the initial observation, ${\bf x_0}$,  towards a final state, ${\bf X^*}$, predicted by a model is expressed in the pink color in the same figure. The uncertainty of the model in representing the observations can be expressed in a number of ways. It may be defined by an absolute error $ E $, which measures distance, in a certain metrics, between the final observation ${\bf x^*}$, which is considered a  target and the  computed prediction ${\bf X^*}$. 
This  work proposes measures  for this error and links the proposed uncertainty quantifier 
with dynamical objects present in the model \eqref{v(t)}. This is done without making any mathematical assumption about the observations.

This paper is structured as follows. In section 2, we will discuss and develop   an  approach to uncertainty quantification recently taken for an ocean study case in \cite{garciasanchez2020}.  We will show evidence of connections between the introduced definitions and 
the dynamical  objects of system \eqref{v(t)}. In section 3, we will  provide formal results to show how, for specific simple examples of the system \eqref{v(t)}, the definition of uncertainty quantification given in  section 2, is able to highlight these dynamical objects. 
These results support one of the main findings of this article which is that stable invariant manifolds of hyperbolic trajectories provide a {\em structure}  to forward uncertainty quantification. We will see in detail what is meant by  this statement.  Section 4 presents a discussion with further examples that illustrate the findings of this work and links   uncertainty quantification to Lagrangian Descriptors and other Lagrangian indicators found in the literature. Finally, in section 5,  we will provide the conclusions.



\begin{figure}[htb!]
  \begin{center}
 \includegraphics[scale=0.5]{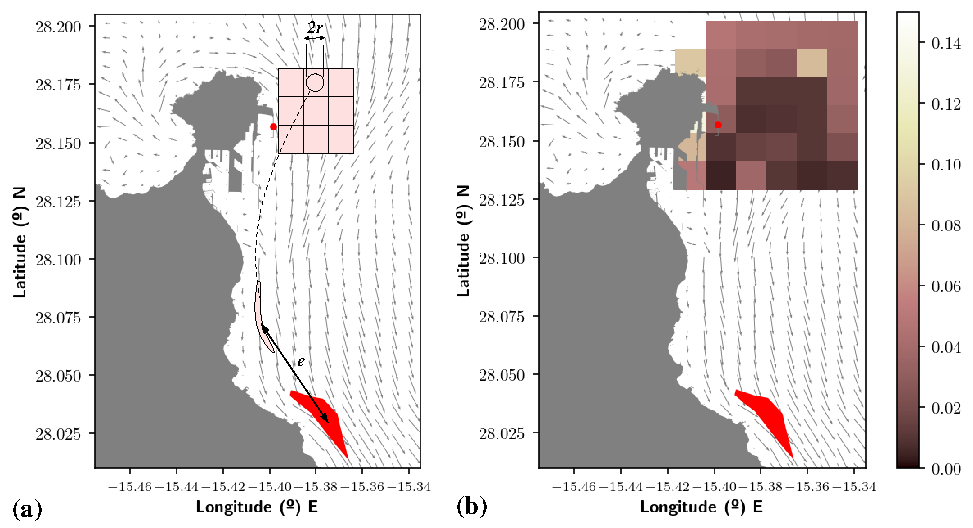} 
  \end{center}
  \caption{A graphical representation of the uncertainty quantification in a neighbourhood of the accident location (red dot) with respect to a target (the red blob). a) A schematic representation of Eq. \eqref{metric} computation for an initial pink blob with radius $r$ and its evolution; b) results of the computation of  Eq. \eqref{metric}  for different initial blobs of radius $r=3.5\cdot 10^{-3}(\degree)$ placed in the center of each cell in the meshgrid. The colormap express in degrees values for the final distances to the target.}
  \label{fig:blobse}
\end{figure}

\section{Transport uncertainty quantification in ocean models}
In this section, we propose uncertainty quantification metrics. We work in the context  of transport in ocean models, starting the discussion from the results presented in \cite{garciasanchez2020},  which considers the performance of very high resolution tools for the monitoring and assessment of environmental hazards in coastal areas. In particular, this work, following the scheme presented in figure \ref{fig:uqsk}, quantifies uncertainties of several high resolution hydrodynamic models in the area of Gran Canaria, by measuring an error between the observations and predictions for the evolution of a diesel fuel spill event, well documented by port authorities and tracked with very high resolution remote sensing products. The pollution event was produced after  the collision  of the passenger ferry `Volc\'an de Tamasite' with the Nelson Mandela dike in La Luz Port on April 2017. After the crash, supply pipes along the dike were broken and diesel fuel poured into the sea. SAR Sentinel 1 images  were available in the area  approximately one day and a half after the accident from which it was possible to identify the spill. For that period ocean currents in the area were available from two sources. One was the Copernicus Marine Service model for the Iberian-Biscay-Irish  region (CMEMS IBI-PHY, IBI hereafter)  and other was a very high resolution model of Puerto de la Luz setup by Puertos del Estado  currently implemented and running operationally in different Spanish Port Authorities within the SAMOA port forecast system \cite{sotillo2020coastal}.  In \cite{garciasanchez2020}, for each  model, uncertainties are quantified by measuring an error with respect to a target, a `ground truth', which should be recovered by the model. In particular, following the scheme proposed in figure \ref{fig:uqsk}, figure \ref{fig:blobse} illustrates how the uncertainty quantifier metrics is proposed.   The ferry impact point, the Nelson Mandela dike, which is marked with the red dot,  represents an initial observation. The currents represented in the background are the ones obtained from the model set by Puertos del Estado at day $t_0$. The elongated  red blob represents the spill as identified from satellite images one day and a half after the accident. This is the target observation, the `ground truth'. 
In panel (a) a  set of initial conditions, ${\bf X_0}$, are selected in the neighbourhood of the initial observation and marked in a rectangular pink domain. This domain is divided into sub-domains and in each one a small circular blob with radius  $r=3.5\cdot 10^{-3} (\degree)$, is  evolved in a time interval from the initial time $t_0$. This is represented just for one of the sub-domains. The blob distorts while  it approaches the target observed spill. A way to measure the `proximity' between the blob and the target spill at the final time, $\textit{t}^*$, is to compute the distance between the centroid of the evolved blob, ($\textbf{c}_m$), and that of the observed spill, ($\textbf{c}_g^*$). Mathematically this can be expressed as:

\begin{equation}
  e(t^*) = \|\textbf{c}_m(t^*) - \textbf{c}_g^*\|.
  \label{metric}
\end{equation}

\noindent
Here, $\|\cdot \|$ denotes the modulus of the vector. The position of the centroid of the evolved spill, $\textbf{c}_m$, depends on time $\textit{t}$, while  the centroid of the ground value slick $\textbf{c}_g^*$ does not, because it is an observation at a final time, $t^*$.
The centroid of a finite set of $N$ points $\lbrace \mathbf{x}_{k}\rbrace_{k\in \mathbb{N}} \in \mathbb{R}^n$ is defined as:

\begin{equation}
    \textbf{c} = \dfrac{1}{N} \sum_{k = 1}^N \mathbf{x}_k \,,
    \label{centroid}
\end{equation}

\noindent
where ${\bf x}_k$ are the (lon, lat) coordinates that define the contour of the slick in an equirectangular projection. The solid black line $e$ in panel (a) is the distance between centroids of the modelled and observed blobs. This is the value taken by Eq. \eqref{metric}.
The procedure is repeated for additional circles centered in different positions of the sub-domains or cells,  in the neighbourhood of the accident place.
The evolution of  each of these circular blobs is different due to the chaotic nature of transport in this setting. Panel (b) in figure \ref{fig:blobse} provides a visualization on how this calculation changes in each sub-domain. This panel represents  a  colormap which is placed in the port neighbourhood, at the same position than the pink domain in panel (a). It  considers a $6\times 6$ mesh-grid, and the colorcode within each cell represents the value of Eq. \eqref{metric}  obtained for blobs with initial positions in the center of each cell and initial radius of $r=3.5\cdot 10^{-3}(\degree)$. We observe a nonuniform distribution.
\begin{figure}[htb!]
  \begin{center}
  a)\includegraphics[scale=0.5]{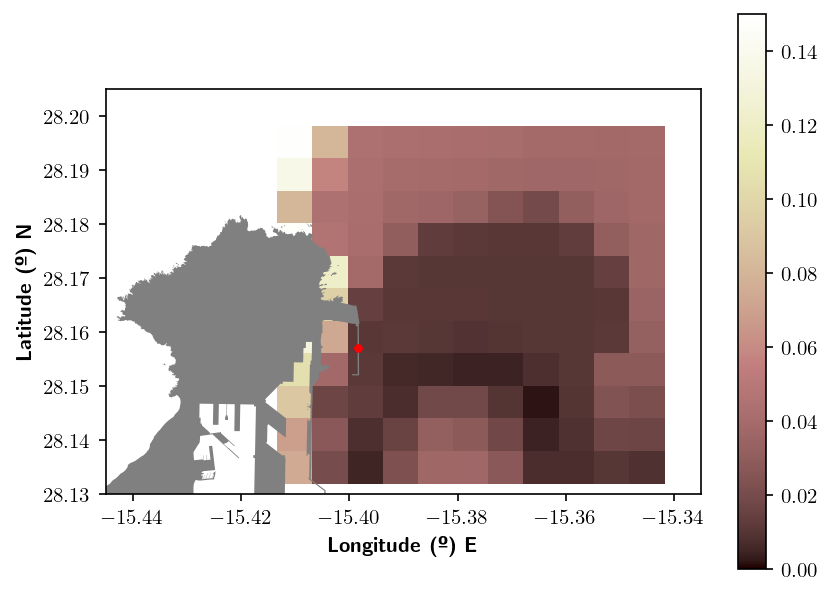}
  b)\includegraphics[scale=0.5]{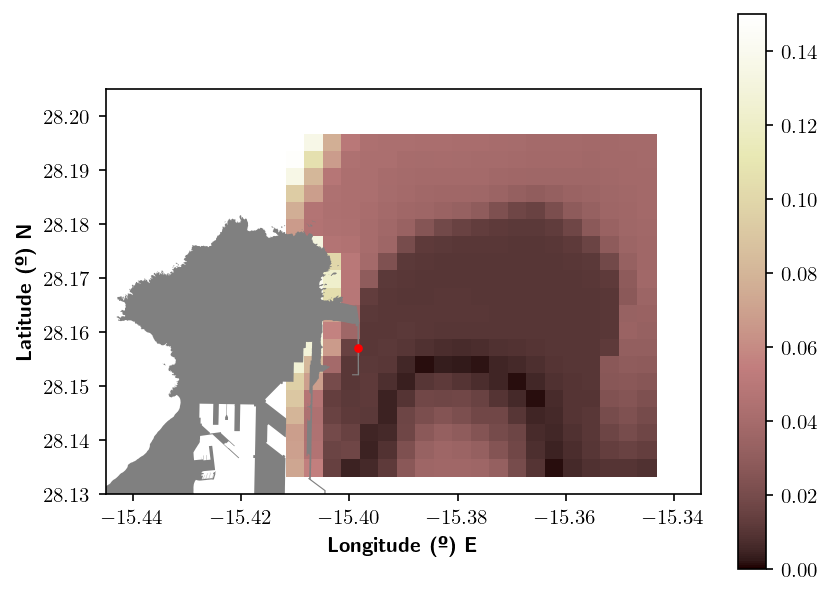}
 \end{center}
  \caption{A quantitative measure of uncertainty for the Puertos del Estado model for la Luz Port, based on errors measured as distances to the target 'ground truth' spill. Higher colour-map values correspond to larger errors or uncertainties. The accident location is marked with a red asterisk. a) Uncertainties computed on a $11\times 11$ mesh considering initial blobs with  radius $r=3.5\cdot 10^{-4}(\degree)$; b) uncertainties associated with a $20\times 20$ mesh considering initial blobs with radius $r=3.5\cdot 10^{-5}(\degree)$.}
  \label{fig:blobs}
\end{figure}

\begin{figure}[htb!]
  \begin{center}
  a)\includegraphics[scale=0.5]{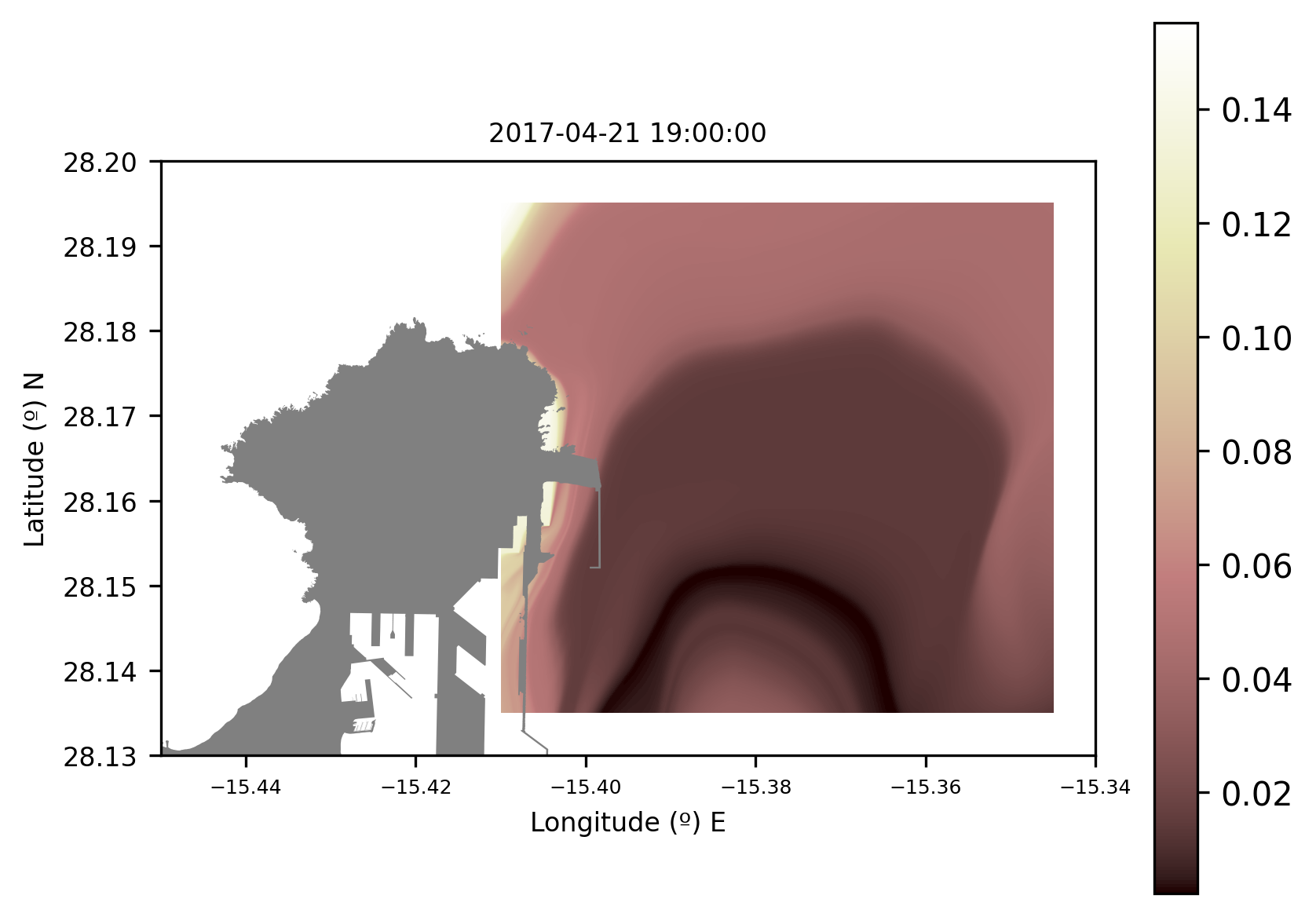} \\
  b)\includegraphics[scale=0.4]{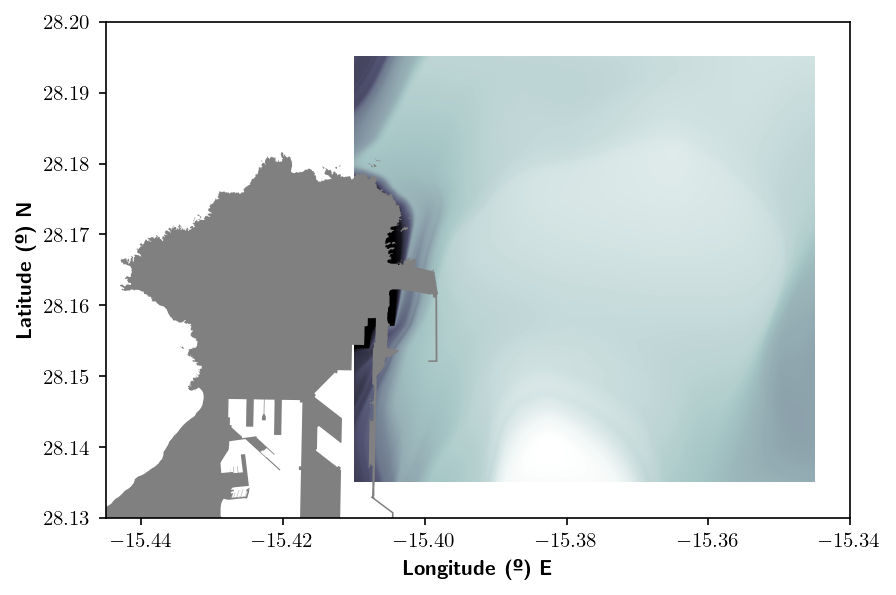}
  c)\includegraphics[scale=0.4]{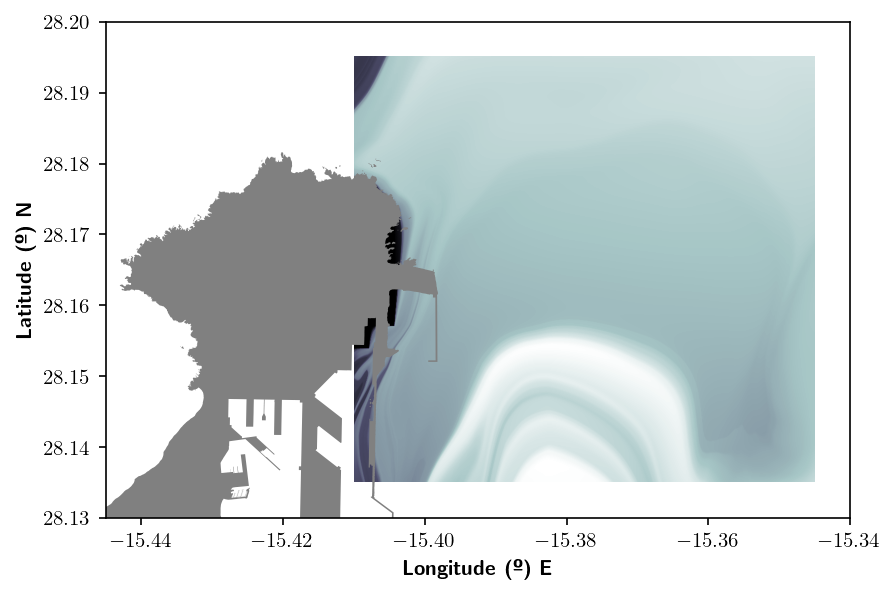}
  \end{center}
  \caption{a) Evaluation of the Lagrangian uncertainty quantification given in equation \eqref{metric2}; b) stable manifolds as revealed by the  $M$ function \eqref{metricM} using $\tau=1.5$ days; c) stable manifolds as revealed by the  $M$ function using $\tau=2$ days  }
  \label{fig:vt manifolds}
\end{figure}
\begin{figure}[htb!]
  \begin{center}
  \includegraphics[scale=0.5]{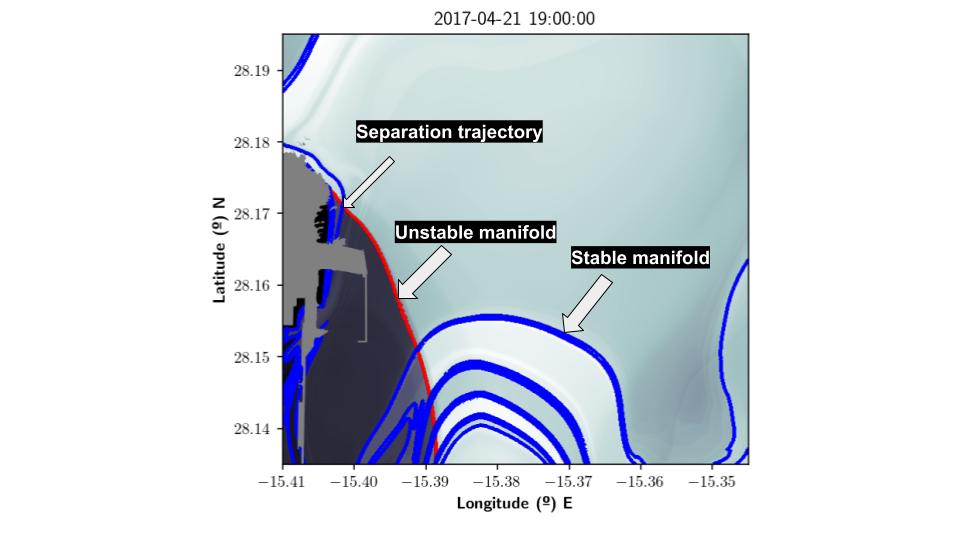} 
  \end{center}
  \caption{ The Lagrangian skeleton of the ocean model on the accident day, the 21 April 2017 at 19:30h. LCSs are highlighted by the red  and blue colors over the  gray in the background.}
  \label{fig:manifolds2}
\end{figure}

 
 Figure \ref{fig:blobs} expands results of Figure \ref{fig:blobse}b) by
increasing the number of cells in the neighborhood domain  and decreasing the  radius of each initial blob. In particular, Figure \ref{fig:blobs} a) considers a mesh  $11\times 11$ and initial blobs with radius $r=3.5\cdot 10^{-4}(\degree)$. Figure \ref{fig:blobs} b)  increases the number of cells to a $20\times 20$ mesh-grid and decreases the radius of initial blobs to $3.5 \cdot 10^{-5}(\degree)$.  Figure \ref{fig:blobs}b) makes visible an underlying structure which is directly related to the invariant dynamical structures, as we will shown later. In the limit $r\to 0$, expression \eqref{metric} is rewritten as:

 \begin{equation}
  L_{UQ}(t^*) = \|\textbf{x}(t^*) - \textbf{c}_g^*\|.
  \label{metric2}
\end{equation}

\noindent
where $\textbf{x}(t)$ is a trajectory of a fluid parcel, a solution to the system  \eqref{v(t)}, with initial position $\textbf{x}_0$ at each cell in the grid. Figure \ref{fig:vt manifolds}a) illustrates the results of this calculation in a very fine grid. 

One of the goals of this paper is to establish connections between the uncertainty quantifier \eqref{metric2} and invariant dynamical objects that control transport in vector fields. 
These  are geometrical objects  that organize particle trajectories schematically into regions corresponding to qualitatively distinct dynamical behaviors. In the context of fluid dynamics these objects are referred to as Lagrangian Coherent Structures (LCS). An  essential  ingredient  of  the LCS  are  hyperbolic  trajectories characterized  by  high  contraction  and  expansion  rates.   Directions  of  contraction  and  expansion  define, respectively, stable and unstable directions, which are, respectively, related to the stable and unstable manifolds.   In the context of the ocean model and the event described above, Garcia-S\'anchez et al. have shown in \cite{garciasanchez2020} that there exists a  hyperbolic trajectory located on the coastline very close to the accident point, in  a detachment configuration, which is related to the phenomena of flow separation.  Under this configuration, the stable manifold of the hyperbolic trajectory is aligned with the coast, and the unstable manifold is transversal to it. Figure \ref{fig:manifolds2} illustrates, respectively, in blue and red the stable and unstable manifolds of the separation trajectory for this particular case that we studied. Any blob placed in the neighbourhood of the separation trajectory,  eventually evolves to become aligned with the unstable manifold, which is an attracting material curve. The configuration of the unstable manifold in Figure \ref{fig:manifolds2}, is consistent with the observed spill marked in red in Figure \ref{fig:blobse}. Indeed, the observed spill has evolved to be completely aligned with the unstable manifold. The invariant dynamical structures displayed in Figure \ref{fig:manifolds2} have been obtained by means  of the Lagrangian Descriptor \cite{madrid2009,mendoza2010,mancho2013}, which measures arc-length of trajectories:

 \begin{equation}
  M= \int_{t_0-\tau}^{t_0+\tau} \left\|\frac{d\textbf{x}(t)}{dt}\right\| dt.
  \label{metricMd}
\end{equation}

\noindent
In \cite{lopesino2017theoretical}, for some generalizations of the expression \eqref{metricMd}, is proved  that stable and unstable manifolds are aligned with singular features of this function. The forward integration 
highlights the stable manifold, while the backwards integration highlights the unstable manifold.  In this work the notion of singular feature is related to an undefined directional derivative in a direction transverse to the manifold curve. In Figure \ref{fig:manifolds2}, red and blue features are placed on the  singular features of $M$.

The metric given in \eqref{metric2} posses similarities with the forward definition of the function $M$:

 \begin{equation}
  M= \int_{t_0}^{t^*=t_0+\tau} \left\|\frac{d\textbf{x}(t)}{dt}\right\| dt.
  \label{metricM}
\end{equation}

\noindent
The analogy between structures obtained from expressions \eqref{metric2} and \eqref{metricM}, is confirmed from figure  \ref{fig:vt manifolds}.  Panel a) displays  expression \eqref{metric2} results, while panels b) and c) display Eq. \eqref{metricM} results for $\tau=1.5$ and 2 days respectively. Figure \ref{fig:vt manifolds} highlights singularities both for the forward uncertainty quantifier and for the forward $M$ function. From this figure is clear that  Uncertainties displayed in panel a)  reach minimum values along the stable manifold.  Indeed,  the stable manifold is an optimal pathway towards the unstable manifold, which in turn is an attracting material curves towards which all fluid parcels evolve.  In this case, as the observed evolution of the spill is aligned with the unstable manifold, it is expected that minimum values of the  uncertainty quantifier correspond to the stable manifold. On the other hand indicators of {\em structural uncertainty} that show inadequacy of the model, would correspond to situations in which these minimum values are not reached along the stable manifolds. This would have been the case if, for instance, the observed spill would have  been found far from the unstable manifold.


\section{Formal results }
\begin{figure}[htb!]
  \begin{center}
  \includegraphics[scale=0.8]{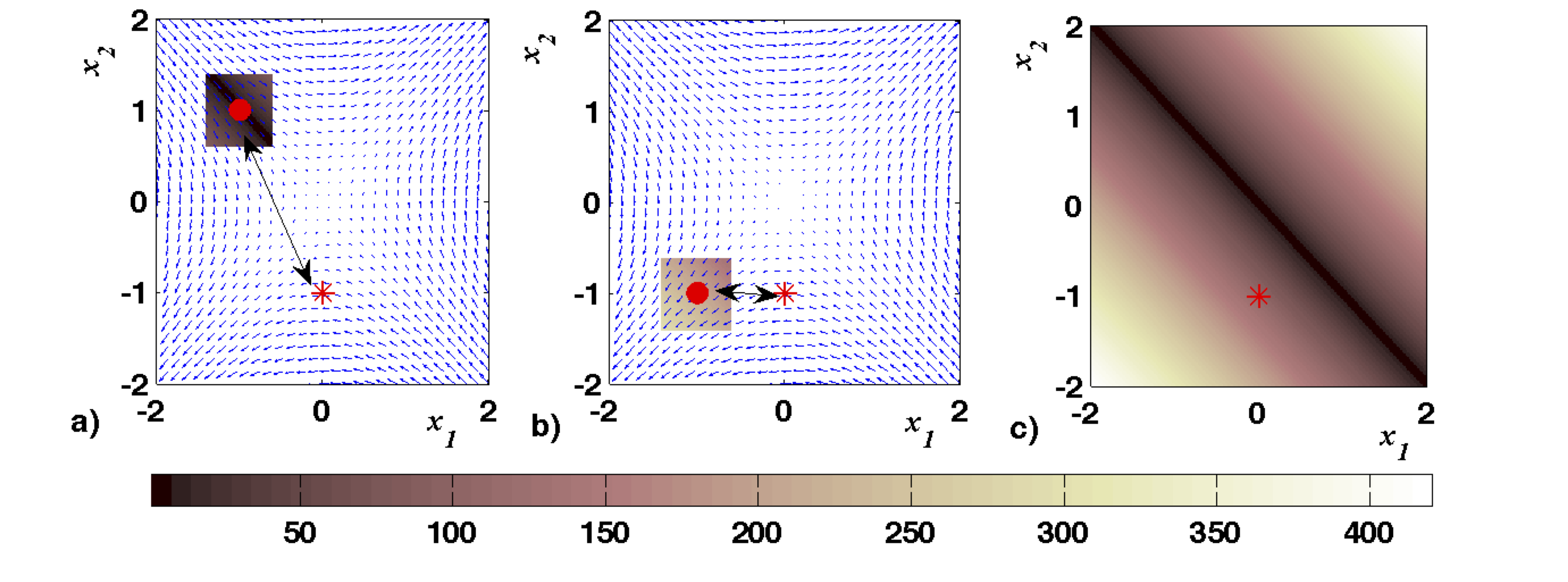} 
  \end{center}
  \caption{a) and b) two successive observations at time $t_0$ (red circle) and time $t^*$ (red asterisk). The transition between them is modeled by the vector field at the background  and $L_{UQ}$ is represented in the neighbourhood of the initial observation. Uncertainties are lower for the case a) as confirmed by the colorbar; c) $L_{UQ}$ evaluated in the whole domain. This expresses the uncertainty associated to the planar flow,  assuming observations that start in different points of the representation domain and that evolve towards  the red asterisk (target position). }
  \label{fig:ex}
\end{figure}

This section is focused on   illustrating the analogies found in figure  \ref{fig:vt manifolds} between the pattern in  panel a), which is  related to uncertainty quantification, and those in panels b) and c), which are related to stable invariant manifolds.   For the event described in the previous section  computations in panel a) are performed with the measure  proposed in Eq. \eqref{metric2} to quantify uncertainty  in a neighbourhood of the accident location. When the calculation is performed on a very fine grid it is found that the uncertainty $L_{UQ}$ has a {\em structure} that is related  to that of invariant stable manifolds. 
In order to find answers to the question of   why there are singular structures obtained from the expression \eqref{metric2} that are aligned with invariant stable manifolds of hyperbolic trajectories, we need to assume an explicit  expression for the vector field in Equation \eqref{v(t)}. Indeed in the example discussed in the previous section, velocities are given as data-sets and this makes difficult to proceed with exact calculations. In this way, in this section  we perform an analysis in simplified planar vector fields with exact explicit expressions. Figure \ref{fig:ex} displays such example. In panel a) an initial observation at time $t_0$ is marked with a red circle. At a later time, $t^*$, this observation is at the target position marked with a red asterisk.  We model this evolution with the vector field at the background, which is represented at time $t_0$.  In order to quantify the uncertainty in the neighbourhood  of the initial observation, we display  the evaluation of     $L_{UQ}$ around it. The colorbar placed at the  bottom of the figure, measures the uncertainty in the domain units. This colorbar indicates that uncertainties in this neighbourhood are low. At least lower than those displayed for the experiment in panel b). In this case the initial observation is at a different position, although the target is kept at the same position 
than in a). In this case uncertainties are larger, indicating that the vector field in the background is not a good model to express the transition between these two observations. Finally panel c)    
displays  $L_{UQ}$ evaluated in the whole domain. This expresses the uncertainty associated to the planar flow,  assuming that observations  start in different points of the representation domain and that they evolve towards the assumed observed target position at the red asterisk. The panel confirms that model has a better performance for those initial observations placed at the diagonal. Results presented next in this section are displayed following this type of representation.  

Our analysis in this section follows the spirit of the work by \cite{mancho2013,lopesino2017theoretical,lopesino2015lagrangian,garcia2018detection}.  We assume the definition of singular features given there, by considering that these are features of  $L_{UQ}$  on which the transversal derivative is not defined. We will prove, for simple selected examples, that stable manifolds are aligned with those singular features of $L_{UQ}$.

Finally, before beginning our discussion the definition \eqref{metric2} is generalized as follows:

\begin{equation}
    L_{UQ}({\bf x}_0, t, p, \delta) = \left[\sum_{i=1}^n\left|x_i(t)-x_{i}^*\right|^p \right]^{\frac{1}{p}}, \ p>1, \delta < 1,  \ {\bf x}_0\in \mathbb{R}^n.
    \label{eq:UQ1}
\end{equation}
Here,  ${\bf x}_0$ is the initial condition of the trajectory $(x_1(t),x_2(t),..,x_n(t))$. In the particular case described in the previous section $n=2$, which  corresponds to the ocean surface, and $p=2$. The coordinates of the target ${\bf c}_g^*$ are $(x_{1}^*,x_2^*)$ An alternative for this expression that we will use is the following:
\begin{equation}
    L_{UQ}({\bf x}_0, t, p) = \sum_{i=1}^n\left|x_i(t)-x_{i}^*\right|^p, \ p\leq 1, \ {\bf x}_0\in \mathbb{R}^n.
    \label{eq:UQ2}
\end{equation}


\subsection{The autonomous saddle point}
The first example that we analyze is the vector field that corresponds to the Hamiltonian linear saddle case for which ${\bf x}\in\mathbb{R}^2 $, satisfying ${\bf x}=(x,y)$ and the equations of motions are:
\begin{eqnarray}
   && \frac{dx}{dt}  = \lambda x, \nonumber\\ 
  &&  \frac{dy}{dt}  = -\lambda y,
    \label{eq:saddle}
\end{eqnarray}
where  $\lambda > 0 $. For any initial observation $(x_0 , y_0 )$, we consider the unique solution of this system passing through the condition $(x_0 , y_0 )$, which is:
\begin{equation}
    \begin{cases} 
    x(t)  = x_0 e^{\lambda t} \\ 
    y(t)  = y_0e^{-\lambda t},
    \end{cases} \quad \lambda > 0 
    \label{eq:saddle solution}
\end{equation}

For this example, the origin $(0,0)$ is a hyperbolic fixed point with stable and unstable manifolds:
\begin{equation}
    W^s(0,0) = \{ (x,y)\in \mathbb{R}^2 \colon x=0, y\neq 0 \},
    \label{manifold: stable saddle}
\end{equation}
\begin{equation}
    W^u(0,0) = \{ (x,y)\in \mathbb{R}^2 \colon x\neq 0, y= 0 \},
    \label{manifold: unstable saddle}
\end{equation}

For simplicity we assume, without loss of generality, that $t_0=0$ (this is possible for autonomous systems). We consider the target in the position $(x^*,y^*)$ and apply \eqref{eq:saddle solution} to \eqref{eq:UQ2} to obtain:
\[
L_{UQ}(\textbf{x}_0, t, p) = |x_0 e^{\lambda t} - x^*|^p + |y_0e^{-\lambda t} - y^*|^p. 
\]

Regrouping terms, we get

\begin{equation}
L_{UQ}(\textbf{x}_0, t, p) = |x_0|^p \omega^{-p} |1 - a \omega|^p + |y_0|^p  |\omega - b |^p, \ \textrm{where} \ a = \frac{x^*}{x_0}, b = \frac{y^*}{y_0}.
\label{saddleContinuous solution}
\end{equation}
In this expression, $\omega=e^{-\lambda t}$, which always satisfies $\omega>0$. We explore separately the first and second terms.
For the factor $ |1 - a \omega|^p$  there exists a $t_L$ such that if $t>t_L$, then $a \omega<<1$ and
$(1-a \omega)>0$. 
This is always the case for $a<0$ and is a plausible assumption for $a>0$, if $a<<\omega^{-1}$. We recall that $a= x^*/x_0$ and that therefore such $t_L$ exists if $x_0 \neq 0$.
In this case positiveness is guaranteed for sufficiently large $t$, i.e, a
 Taylor series around $\omega=0$, attained if $t\gg1$ and $t>t_L$, is performed for the binomial:
\[
(1-a\omega)^p = 1-a p \omega+\frac{1}{2} a^2 (p-1) p \omega^{2}-\frac{1}{6} \omega^{3} \left(a^3 (p-2) (p-1) p\right)+O\left(\omega^{4}\right).
\]
Therefore, 
\[
\frac{1}{\omega^{p}} (1-a\omega)^p 
 = \omega^{-p} + a p \omega^{(1-p)} + \frac{1}{2} a^2 (p-1) p \omega^{(2-p)} + O\left(\omega^{(3-p)}\right).
\]
We recall that $p\leq 1$ and $\omega=e^{-\lambda t}$. This yields,
\[
    |x_0|^p|e^{\lambda t}-a|^p =  |x_0|^p e^{\lambda t p} + O\left(|x_0|^p a e^{-\lambda t(1-p)}\right).
\]

We analyse next the second term in Eq.\eqref{saddleContinuous solution}. The sign of $(\omega-b)$ depends crucially on the sign of $b$ given that for $t\gg 1$ the term $\omega $  can be as small as we like. Let us consider $t_L$ such that if $t>t_L$  and $b<0$, then  $(e^{-\lambda t}-b)>0$. The following is valid for $t$ above this lower value:
\[
|\omega-b| ^p =(\omega -b)^p = (-b)^p + p \omega (-b)^{p-1} + \frac{1}{2} (p-1) p \omega^2(-b)^{p-2}
\]
\[
+ \frac{1}{6} (p-2) (p-1) p \omega^3(-b)^{p-3}+O\left(\omega^4\right)
\]
Therefore, when $t\gg 1$, $b<0$ and $e^{-\lambda t} \ll 1$,
\begin{equation}
    (\omega-b)^p = (-b)^p + O\left(e^{-\lambda t}\right), e^{-\lambda t} \ll 1
    \label{dominant positive unstable}
\end{equation}

Now, we consider $b>0$ and $t>t_L$ satisfying $(\omega-b)<0$. Therefore:
\[
|\omega-b|^p=(-\omega+b)^p = b^p + p \omega b^{p-1}+\frac{1}{2} (p-1) p \omega^2 b^{p-2}
\]
\[
+ \frac{1}{6} (p-2) (p-1) p \omega^3 b^{p-3}+O\left(\omega^4\right).
\]
Therefore, when  $t\gg 1$, $b>0$ and $e^{-\lambda t} \ll 1$,
\begin{equation}
    (-\omega+b)^p = b^p + O\left(\omega\right).
    \label{dominant negative unstable}
\end{equation}
Finally, we conclude that
\begin{equation}
    |\omega-b|^p = |b|^p+  O\left(e^{-\lambda t}\right), \, \omega \ll 1.
    \label{dominant abs unstable}
\end{equation}
Therefore, since $b = y^*/y_0$
\[
    |y_0|^p|e^{-\lambda t}-b|^p =  |y^*|^p+  O\left(|y_0|^p e^{-\lambda t}\right), \, e^{-\lambda t} \ll 1.
\]
Thus, we can approximate $L_{UQ}$ as
\begin{equation}
L_{UQ}(\textbf{x}_0, t, p) \approx |x_0e^{\lambda t}|^{p} + |y^*|^p = |x_0|e^{\lambda t p} + |y^*|^p \label{luqapp1} 
\end{equation}
where the dominant term is $|x_0|^p e^{\lambda tp}$. Hence, to leading order, the stable manifold at $x=0$ is aligned with a singular feature of $L_{UQ}$  for `sufficiently large'  $t$. This statement however must be considered in the sense that expression \label{luqapp1} expresses a good approximation for $L_{UQ}$ as far as $x_0\neq 0$, and that for any $|x_0|> 0$  is valid for a sufficient large $t$ satisfying, $t>t_L$.



\begin{figure}[htb!]
  \begin{center}
  a)\includegraphics[scale=0.5]{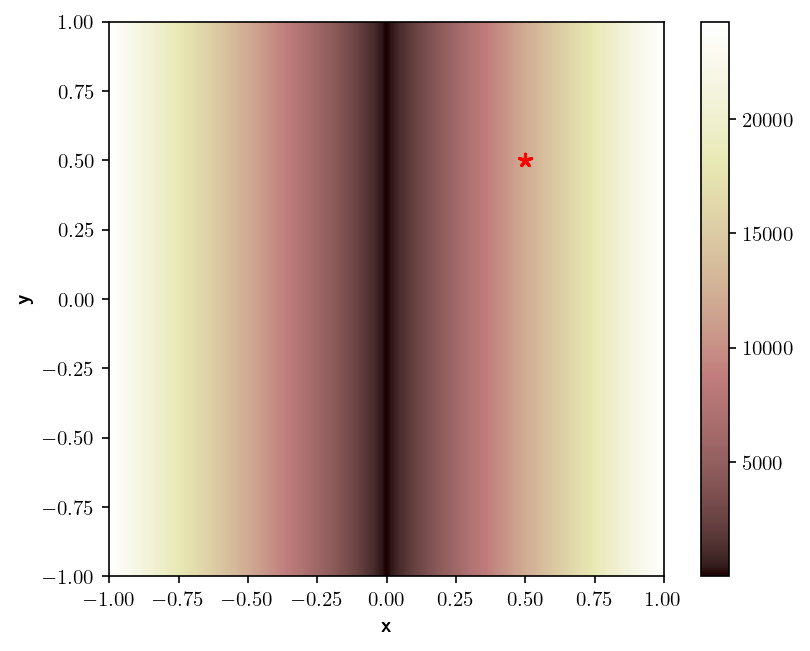} 
  b)\includegraphics[scale=0.5]{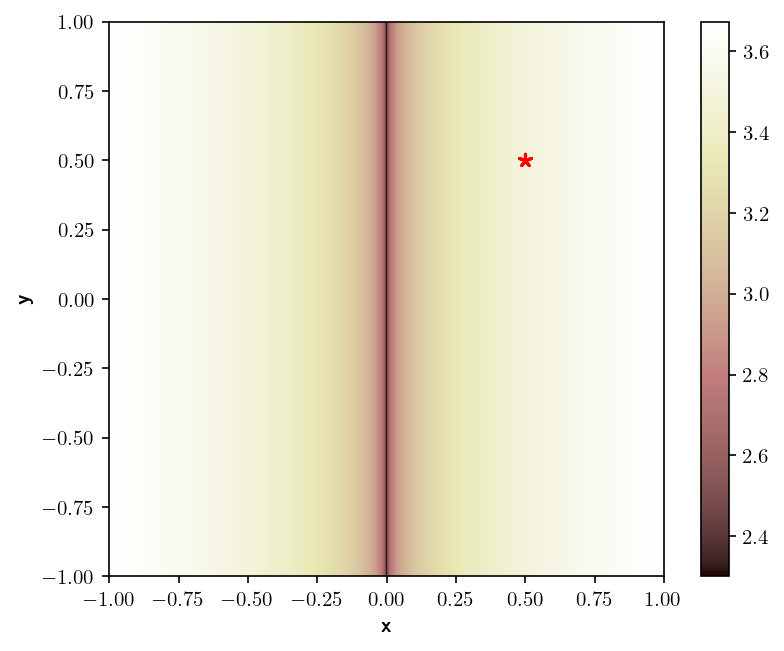}
  \end{center}
  \caption{A representation of Eq.(\ref{eq:UQ1}) and Eq.\eqref{eq:UQ2}  for $t^*=10$ and target $(x^* , y^* )= (0.5,0.5)$ for the Hamiltonian linear saddle vector field. a) $ p=2$; b)  $p=0.1$ . It can be appreciated how the stable manifold is aligned with a singular feature. }
  \label{fig:continuous saddle LUQ}
\end{figure}

We have shown that Eq. \eqref{eq:UQ2} is able to highlight the stable manifolds for the autonomous saddle point. Next, we analyse what happens when we consider Eq.\eqref{eq:UQ1}, and in which cases it provides information about the stable manifold of the autonomous saddle point. When we apply \eqref{eq:saddle solution} to \eqref{eq:UQ1}, it yields
\[
L_{UQ}(\textbf{x}_0, t, p) = \left[|x_0 e^{\lambda t} - x^*|^p + |y_0e^{-\lambda t} - y^*|^p\right]^{1/p}.
\]
We remark that our following calculations will be restricted to $p>1$ and integer. Rewriting the expression for $\omega=e^{-\lambda t}$
\begin{eqnarray}
L_{UQ}(\textbf{x}_0, t, p) = \left[\left|\frac{x_0}{\omega}\right|^p |1 - a\omega|^p + |y_0|^p |\omega - b|^p\right]^{1/p}
= \left|\frac{x_0}{\omega}\right| \left[ |1 - a\omega|^p + 
\left| \frac{y_0 }{x_0} \right|^p  \omega^p |\omega - b|^p\right]^{1/p}
\label{eqLUQ2}
\end{eqnarray}
Recalling  that $p>1$, $\omega>0$  and that $(1 - a\omega)>0$ for sufficiently small $\omega$, i.e. sufficiently large $t$,  and the case $b<0$ and $(\omega - b)>0$, a Taylor series around $\omega=0$, for $p$ integer number, for the second factor: 
\begin{eqnarray}
L_{UQ}(\textbf{x}_0, t, p) = \left|\frac{x_0}{\omega}\right| \left[   1- a \omega +\left| \frac{y_0 }{x_0} \right|^p (-b)^p \frac{\omega^p}{p} +O\left(\omega\right)^{p+1} \right]
\end{eqnarray}
Alternatively, considering  the case $b>0$ and $(\omega - b)<0$, a Taylor series around $\omega=0$ for the second factor is as follows: 
\begin{eqnarray}
L_{UQ}(\textbf{x}_0, t, p) = \left|\frac{x_0}{\omega}\right| \left[   1- a \omega +\left| \frac{y_0 }{x_0} \right|^p (b)^p \frac{\omega^p}{p} +O\left(\omega\right)^{p+1} \right]
\end{eqnarray}
Therefore in general:
\begin{eqnarray}
&&L_{UQ}(\textbf{x}_0, t, p) = \left|\frac{x_0}{\omega}\right|
\left[   1- a \omega +\left| \frac{y_0 }{x_0} \right|^p |b|^p \frac{\omega^p}{p} +O\left(\omega\right)^{p+1} \right]
\nonumber\\
&&=  \left|\frac{x_0}{\omega}\right| -a |x_0|+ \left|y_0 \right|^p \left|x_0 \right|^{1-p} \frac{|b|^p}{p} \omega^{p-1} +O\left(\omega\right)^{p} \label{sing}
\end{eqnarray}
We notice that since $p>1$, the terms $\omega^{p-1},\omega^{p} $ have positive exponents and if $\omega \ll 1$ in principle are much smaller than the first term  $\omega^{-1}$ and therefore:
\begin{eqnarray}
L_{UQ}(\textbf{x}_0, t, p) \sim     \frac{1}{p} \left|\frac{x_0}{\omega}\right|, \,\, \omega \ll 1  \label{sing2}
\end{eqnarray}
Hence, the stable manifold  at $x=0$ is aligned with a singular feature of this approximation to  $L_{UQ}$ for `sufficiently large' $t$. This statement though, must be considered with care, because in Eq. \eqref{sing} the terms in $\omega^{p-1},\omega^{p},  $ are multiplied by 
$|x_0|^{1-p}$ which has a negative exponent, and therefore have a singularity at $x_0=0$. For this reason neglecting these terms versus the first one if we are very close to $x_0=0$ would require checking that the products, such as, $|x_0|^{1-p} \omega^{p-1}$ are really small. This implies that  Eq. \eqref{sing2} is correctly approximating $L_{UQ}$ for $\omega\ll 1$
as far as we are sufficiently away from 0 in $x_0$, i.e. $|x_0|^{1-p} \omega^{p-1}\ll 1 \implies \omega^{p-1} \ll |x_0|^{p-1} $. In practice this condition is satisfied for any grid $(x_0, y_0)$ used in later figures that exclude $x_0=0$.  Also we can state  that for sufficiently large $t$, $L_{UQ}$ is very close in almost all the domain to the function given  in Eq. \eqref{sing2}, which possess a "singular feature" aligned with the stable manifold. 
Fig.\ref{fig:continuous saddle LUQ} a) and b),  illustrates how  the stable manifold is aligned with  singular features of the Lagrangian uncertainty quantifier defined either by Eq.(\ref{eq:UQ1}) or (\ref{eq:UQ2}). 

\subsection{The autonomous rotated saddle point}
This second case that we explore is  the vector field of the rotated linear saddle for which the equations of motion are: 
\begin{equation}
    \begin{cases} 
    \dot{x}  = \lambda y \\ 
    \dot{y}  = \lambda x,
    \end{cases} \quad \lambda > 0
    \label{eq:rotated saddle}
\end{equation}
The general solution to this system is:
\begin{equation}
    \begin{cases} 
    x(t)  = a e^{\lambda t} + b e^{-\lambda t} \\ 
    y(t)  = a e^{\lambda t} - b e^{-\lambda t},
    \end{cases} \quad \lambda > 0
    \label{eq:rotated saddle solution}
\end{equation}
where $a$ and $b$ depend on the initial conditions $x_0$ and $y_0$ as follows: 
\[
a = \frac{x_0 + y_0}{2}, \ b = \frac{x_0 - y_0}{2}.
\]
Here, $a = 0$ corresponds to the stable manifold of the hyperbolic fixed point placed at the origin and $b = 0$ corresponds to its unstable manifold. 

For any initial observation $(x_0, y_0)$ and final target  observation $(x^*, y^*)$,  we introduce the solution \eqref{eq:rotated saddle solution} into Eq.\eqref{eq:UQ2} obtaining: 
\begin{equation}
    L_{UQ}(\textbf{x}, t, p) = |a e^{\lambda t}  + b e^{-\lambda t} - x^*|^p + |a e^{\lambda t} -  b e^{-\lambda t} - y^*|^p
    \label{saddleRotated continuous}
\end{equation}

We expand next the first term in Eq.\eqref{saddleRotated continuous}. In particular, we consider the case in which $a>0$ for which it is always possible to find a $t>t_L$ in which $(a e^{\lambda t}  + b e^{-\lambda t} - x^*)>0$:
\[
|a e^{\lambda t}  + b e^{-\lambda t} - x^*| ^p = (a e^{\lambda t}  + b e^{-\lambda t} - x^*)^p = \left(\frac{a}{\omega} + b \omega - x^*\right)^p 
\]
\[
 = \omega^{-p}\left(a + b\omega^{2} - x^*\omega\right)^p,
\]
Here we have used $\omega=e^{-\lambda t} $. The Taylor expansion of the binomial when $\omega \ll 1$, i.e., when $t\gg 1$, is

\[
\omega^{-p}\left(a + b\omega^{2} - x^*\omega\right)^p =\omega^{-p}\left[ a^p- p\omega x^* a^{p-1}+\omega^{2} \left(b p a^{p-1}+\frac{1}{2} (p-1) p (x^*)^2 a^{p-2}\right)\right.
\]
\[
    \left. -\frac{1}{6} \omega^{3} \left((p-1) p x^* a^{p-3} \left(6 a b+p (x^*)^2-2(x^*)^2\right)\right)+O\left(\omega^{4}\right)\right].
\]

Therefore,
\[
\omega^{-p}\left(a + b\omega^{2} - x^*\omega\right)^p = a\omega^{-p} - p\omega x^* a^{p-1}+\omega^{(2-p)} \left(b p a^{p-1}+\frac{1}{2} (p-1) p (x^*)^2 a^{p-2}\right)
\]
\[
 -\frac{1}{6} \omega^{(3-p)} \left((p-1) p x^* a^{p-3} \left(6 a b+p (x^*)^2-2(x^*)^2\right)\right)+O\left(\omega^{(4-p)}\right)
\]
Since $p<1$, $\omega^{(n-p)}\ll 1$ for $n>1$  when $t\gg 1$. Hence,
\[
    \omega^{-pt}\left(a + b\omega^{2t} - x^*\omega^t\right)^p = a^p\omega^{-pt} + O\left(\omega^{(1-p)t}\right).
\]
By removing the change of variable, we  obtain
\begin{equation}
    \left(a e^{\lambda t}  + b e^{-\lambda t} - x^*\right)^p = a^p e^{\lambda pt} + O\left(e^{-(1-p)\lambda t}\right) \approx a^p e^{\lambda pt}, \ \textrm{when}, \ t\gg 1.
\end{equation}

We consider now   the case $a<0$ that implies $\left(a e^{\lambda t}  + b e^{-\lambda t} - x^*\right)<0$ for some sufficiently large $t$. The Taylor series is the same, except for a minus sign. It yields,
\begin{equation}
    \left(a e^{\lambda t}  + b e^{-\lambda t} - x^*\right)^p = -a^p e^{\lambda pt} - O\left(e^{-(1-p)\lambda t}\right) \approx -a^p e^{\lambda pt}, \ \textrm{when}, \ t\gg 1.
\end{equation}
Finally,
\begin{equation}
    |a e^{\lambda t}  + b e^{-\lambda t} - x^*|^p = |a|^p e^{\lambda pt} + O\left(e^{-(1-p)\lambda t}\right) \approx |a|^p e^{\lambda pt}, \ \textrm{when}, \ t\gg 1.
\end{equation}

We analyse next the second term in Eq.\eqref{saddleRotated continuous}. As before, we start considering the case $a>0$ and $(a e^{\lambda t}  - b e^{-\lambda t} + y^*)>0$. After making the change of variable $\omega = e^{-\lambda t}$ we obtain:

\[
\omega^{-p}\left(a - b\omega^{2} - y^*\omega^t\right)^p =\omega^{-p}\left[ a^p-p \omega^t y^* a^{p-1}-\frac{1}{2} \omega^{2} \left(p a^{p-2} \left(2 a b-p(y^*)^2+(y^*)^2\right)\right)\right.
\]
\[
    \left. +\frac{1}{6} (p-1) p \omega^{3} y^* a^{p-3} \left(6 a b-p (y^*)^2+2 (y^*)^2\right)+O\left(\omega^{4}\right)\right].
\]
Since $p<1$, $\omega^{(n-p)} \ll 1$ for $n>1$ when $t\gg 1$. Hence,
\[
    \omega^{-p}\left(a - b\omega^{2} + y^*\omega\right)^p = a^p\omega^{-p} + O\left(\omega^{(1-p)}\right).
\]
By removing the change of variable, we  obtain:
\begin{equation}
    \left(a e^{\lambda t}  - b e^{-\lambda t} - y^*\right)^p = a^p e^{\lambda pt} + O\left(e^{-(1-p)\lambda t}\right) \approx a^p e^{\lambda pt}, \ \textrm{when}, \ t\gg 1.
\end{equation}
We compute the same term but in the case $a<0$ $\left(a e^{\lambda t}  - b e^{-\lambda t} - y^*\right)<0$. The Taylor series is the same, except for a minus sign. It yields,
\begin{equation}
    \left(a e^{\lambda t}  - b e^{-\lambda t} - y^*\right)^p = -a^p e^{\lambda pt} - O\left(e^{-(1-p)\lambda t}\right) \approx -a^p e^{\lambda pt}, \ \textrm{when}, \ t\gg 1.
\end{equation}
Finally,
\begin{equation}
    |a e^{\lambda t}  - b e^{-\lambda t} - y^*|^p = |a|^p e^{\lambda pt} + O\left(e^{-(1-p)\lambda t}\right) \approx |a|^p e^{\lambda pt}, \ \textrm{when}, \ t\gg 1.
\end{equation}

Thus, 
\begin{equation}
     L_{UQ}=|a e^{\lambda t}  +  b e^{-\lambda t} - x^*|^p + |a e^{\lambda t}  - b e^{-\lambda t} + y^*|^p \approx |a|^p e^{\lambda pt}, \ \textrm{when}, \ t\gg 1. \label{luqrc}
\end{equation}
The stable manifold at $a=0$ is aligned with a singular feature also at $a=0$. 

\begin{figure}[htb!]
  \begin{center}
    a)\includegraphics[scale=0.5]{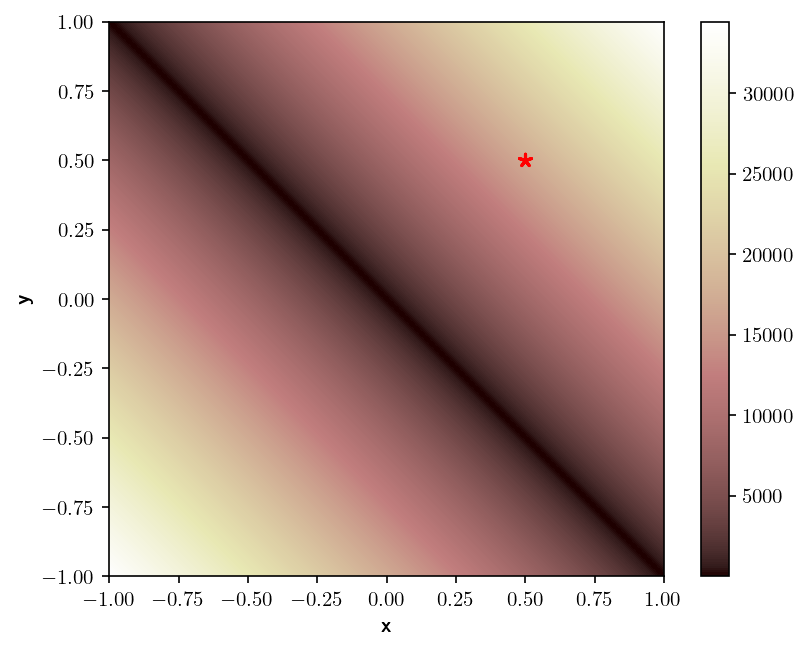}
  b)\includegraphics[scale=0.5]{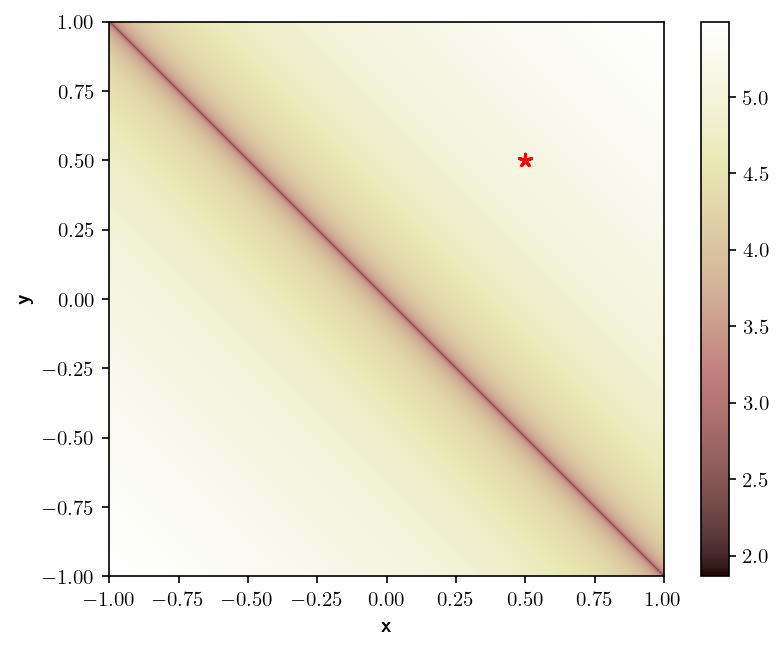} 
  \end{center}
  \caption{A representation of Eq.\eqref{eq:UQ1} and Eq.\eqref{eq:UQ2} for $t^{*} = 10$ and target $(x^{*},y^{*} ) = (0.5,0.5)$ for the rotated linear saddle vector field   for a) $p=2$; b) $p=0.1$.  It can be appreciated how the stable manifold is aligned with a  singular feature. }
  \label{fig:continuous saddle rotated LUQ}
\end{figure}

As for Eq.\eqref{eq:UQ1} considering $p>1$, $p=1/\delta$ and $\omega=e^{-\lambda t}$ leads:
\begin{eqnarray}
   L_{UQ}&=& \left[|a e^{\lambda t}  +  b e^{-\lambda t} - x^*|^p + |a e^{\lambda t}  - b e^{-\lambda t} + y^*|^p\right]^{1/p}\nonumber \\ &=&\omega^{-1}\left[
   |a  +  b \omega^2 - x^*\omega|^p + |a  - b \omega^2 + y^*   \omega|^p
   \right]^{1/p} \nonumber\\
  &=& 2^\frac{1}{p}|a|\omega^{-1}-2^{-1+\frac{1}{p}} (2b+x^*-y^*)+2^{-3+\frac{1}{p}}|a|^{-1} \cdot \nonumber\\
  &&(4 a b + (-1 + p) (-2 b + x^* + y^*)^2) \omega+O(\omega)^2
\end{eqnarray}
As before the leading term in $\omega$ leads to:
\begin{eqnarray}
   L_{UQ}  &\sim& 2^\frac{1}{p}|a|\omega^{-1}
   \end{eqnarray}
and similar considerations apply in neglected terms to the ones made for Eq. \eqref{sing}, but now regarding to singularities at $a=0$.
Fig.\ref{fig:continuous saddle rotated LUQ} a) and b),  illustrate how  the stable manifold is aligned with  singular features of the Lagrangian uncertainty quantifier defined either by Eq.(\ref{eq:UQ1}) or (\ref{eq:UQ2}).

\subsection{Discrete maps}

These findings in the previous two subsections can be  easily extended to discrete time dynamical systems, which are also useful in applications. Discrete time dynamical systems are defined as maps.

\subsubsection*{The autonomous saddle point}

Consider the following linear, area-preserving autonomous map: 
\begin{equation}
    \begin{cases} 
    x_{n+1}  = \lambda x_{n} \\ 
    y_{n+1}  = \frac{1}{\lambda} y_{n},
    \end{cases} \quad \lambda > 1.
    \label{eqmap:saddle}
\end{equation}
For an initial condition $(x_0 , y_0)$, the unique solution of this system is:
\begin{equation}
    \begin{cases} 
    x_n  = x_0 \lambda^n \\ 
    y_n  = y_0 \lambda^{-n},
    \end{cases} \quad \lambda > 1.
    \label{eqmap:saddle solution}
\end{equation}
As for the continuous time case, the origin $(0,0)$ is a hyperbolic fixed point with stable and unstable manifolds:
\begin{equation}
    W^s(0,0) = \{ (x,y)\in \mathbb{R}^2 \colon x=0, y\neq 0 \},
    \label{manifold: stable saddle}
\end{equation}
\begin{equation}
    W^u(0,0) = \{ (x,y)\in \mathbb{R}^2 \colon x\neq 0, y= 0 \},
    \label{manifold: unstable saddle}
\end{equation}
We apply (\ref{eqmap:saddle solution}) to (\ref{eq:UQ2})  to obtain:
\begin{equation}
L_{UQ}(\textbf{x}, n, p) = |x_0 \lambda^n - x^*|^p + |y_0 \lambda^{-n} - y^*|^p. \label{luqmap}
\end{equation}
Regrouping terms, we get
\begin{equation}
L_{UQ}(\textbf{x}, n, p) = |x_0|^p|\lambda^n - a|^p + |y_0|^p|\lambda^{-n} - b|^p, \ \textrm{where} \ a = \frac{x^*}{x_0}, b = \frac{y^*}{y_0}.
\label{saddleDiscrete solution2}
\end{equation}

Considering that the transformation $e^{\lambda t}\to \lambda^n$  can be directly applied to all the calculations performed in the continuous time  case, we recover from Eq.\eqref{saddleContinuous solution}:
\begin{equation}
 L_{UQ} \approx |x_0|^{p}|\lambda^{np}| + O\left(|\lambda|^{-n(1-p)}\right) , \ \textrm{when}, \ n\gg 1.
\end{equation}
When $p < 1$, the dominant term is $|x_0|^{p}\lambda^{np}$. Hence, the stable manifold at $x=0$  is aligned with singular features of  $L_{UQ}$  for `sufficiently large' iteration $n$.

The same is applicable to $L_{UQ}$ \eqref{eq:UQ1}. It yields,
\begin{equation}
    L_{UQ} \approx |x_0||\lambda^{n}| + O\left(1\right) , \ \textrm{when}, \ n\gg 1 \ \textrm{and}\ p\delta < 1.
\end{equation}
We note that the stable manifold at $x=0$ is aligned   with singular features of  $L_{UQ}$ for a `sufficiently large' iteration $n$. The same issues as before regarding singularities on the manifold position apply.

\begin{figure}[htb!]
  \begin{center}
  a)\includegraphics[scale=0.5]{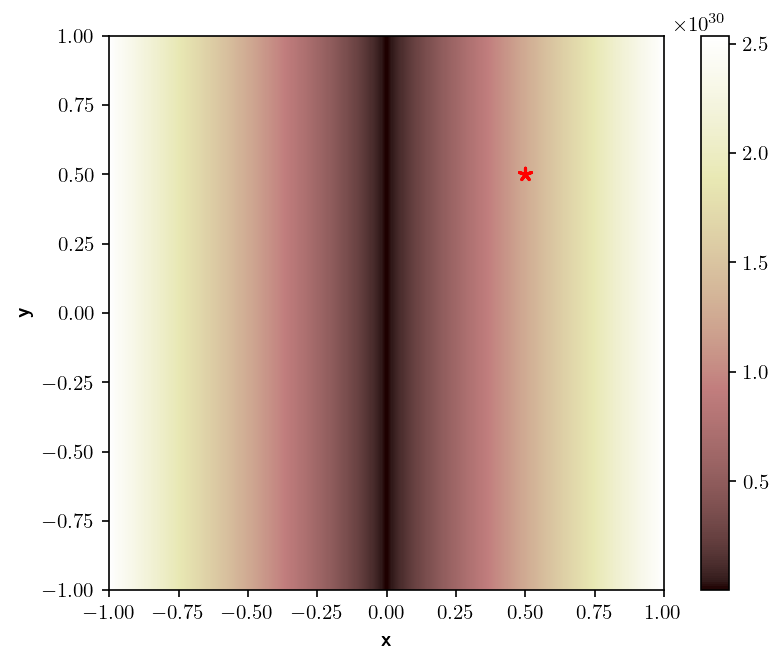}
  b)\includegraphics[scale=0.5]{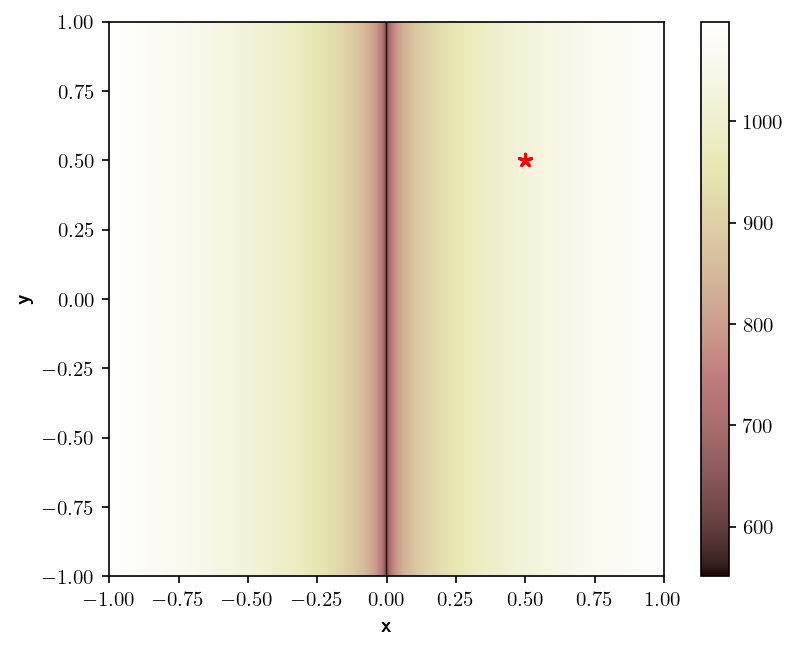} 
  \end{center}
  \caption{
  In Figure a) there is a representation for the autonomous saddle when we apply Eq.(\ref{eq:UQ1}) for $p=2$. In b), we illustrate the same representation but for Eq.(\ref{eq:UQ2}) when  $p=0.1$. It can be appreciated how the stable manifold is aligned with the singular feature. }
  \label{fig:discrete saddle LUQ}
\end{figure}

Fig.\ref{fig:discrete saddle LUQ} a) and b),  illustrate how  the stable manifold is aligned with  singular features of the Lagrangian uncertainty quantifier defined either by Eq.(\ref{eq:UQ1}) or (\ref{eq:UQ2}). 

\subsection*{The autonomous rotated saddle point}

We consider the following discrete dynamical system:
\begin{align}
    F(x,y) &= A\begin{pmatrix}
           x \\
           y
           \end{pmatrix},\  \textrm{where} \  A =
           \frac{1}{2\lambda } 
           \begin{pmatrix}
           \lambda^2 + 1 & \lambda^2 -1\\
           \lambda^2 -1& \lambda^2+1
           \end{pmatrix}, \ \lambda > 1.
          \label{eqmap:rotatedsaddle}
\end{align}
It is easy to see that the stable and the unstable manifolds are given by the vectors $(1,-1)$ and $(1,1)$ respectively. The solution of this system yields to,

\begin{equation}
    \begin{cases} 
    x_n  = a\lambda^{n} + b\lambda^{-n} \\ 
    y_n  = a\lambda^{n} - b\lambda^{-n},
    \end{cases} \quad \lambda > 1
    \label{eqmap:rotated saddle solution}
\end{equation}
where 
\[
a = \frac{x_0 + y_0}{2}, \ b = \frac{x_0 - y_0}{2}.
\]
Again considering the transformation $e^{\lambda t}\to \lambda^n$ and the use of previous results for the continuous time case, we recover from Eq.\eqref{luqrc}:
\[
L_{UQ}(\textbf{x}, n, p) = |a\lambda^n + b\lambda^{-n} - x^*|^p + |a\lambda^{n} - b\lambda^{-n} - y^*|^p.
\]
Therefore, the dominant term is
\[
L_{UQ}(\textbf{x}, n, p) \sim |a|^p|\lambda|^{np}.
\]
Since $a =\frac{x_0 + y_0}{2} $, there is a singular feature at $x=-y$, i.e., in the subspace generated by $(1,-1)$. Hence, the stable manifold is aligned with a singular feature of $L_{UQ}$. 

\begin{figure}[htb!]
  \begin{center}  
  a)\includegraphics[scale=0.5]{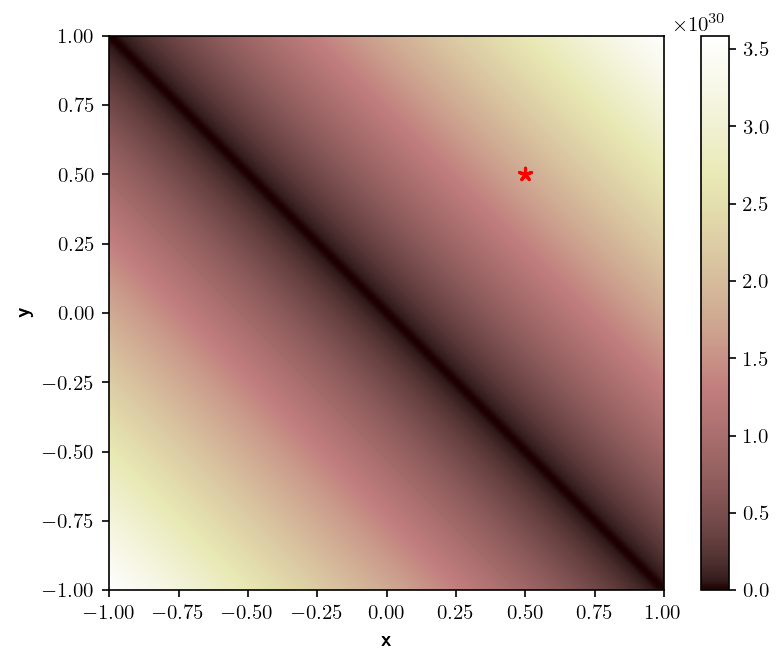}
  b)\includegraphics[scale=0.5]{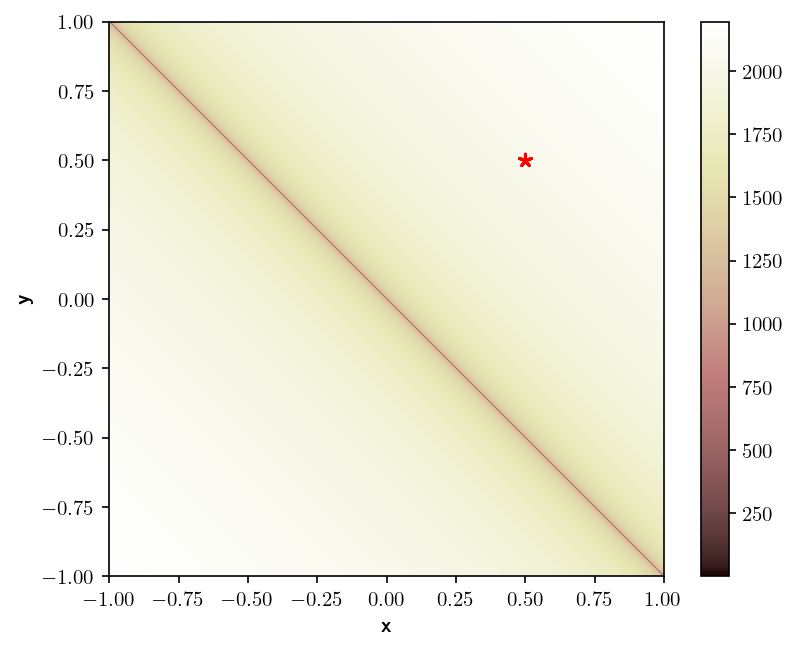} 
  \end{center}
  \caption{In Figure a) there is a representation for the rotated saddle when we apply Eq.(\ref{eq:UQ1}) for  $p=1/ \delta=2$. In b), we illustrate the same representation but for Eq.(\ref{eq:UQ1}) when  $p=1/ \delta=0.1$. It can be appreciate how the singular feature is aligned with stable manifold. }
  \label{fig:discrete saddle rotated LUQ}
\end{figure}

The same is applicable to $L_{UQ}$ \eqref{eq:UQ1}. It yields,
\[
L_{UQ}(\textbf{x}, n, p) \sim |a||\lambda|^{n},
\]
Hence, the stable manifold is aligned with a singular feature of $L_{UQ}$. Fig.\ref{fig:discrete saddle rotated LUQ} a) and b),  illustrates how  the stable manifold is aligned with  singular features of the Lagrangian uncertainty quantifier defined either by Eq.(\ref{eq:UQ1}) or (\ref{eq:UQ2}). 


\section{The Duffing equation}
Results in the previous sections are generalized to the autonomous nonlinear case by means of the Moser's theorem \cite{moser56}. This theorem applies to analytic two-dimensional
symplectic maps having a hyperbolic
fixed point or, similarly, to two-dimensional time-periodic
Hamiltonian vector fields having a hyperbolic
periodic orbit (which can be reduced to the
former case considering a Poincar\'e map). The case of a Hamiltonian nonlinear autonomous system is a one-parameter family of symplectic maps, and therefore Moser’s
theorem applies.  Following proofs sketched by \cite{lopesino2015lagrangian, lopesino2017theoretical} results may be extended to the case of  non-autonomous nonlinear 
dynamical systems   by utilizing results like the Hartman–Gro$\beta$man theorem. 

\begin{figure}[htb!]
  \begin{center}
 a) \includegraphics[scale=0.5]{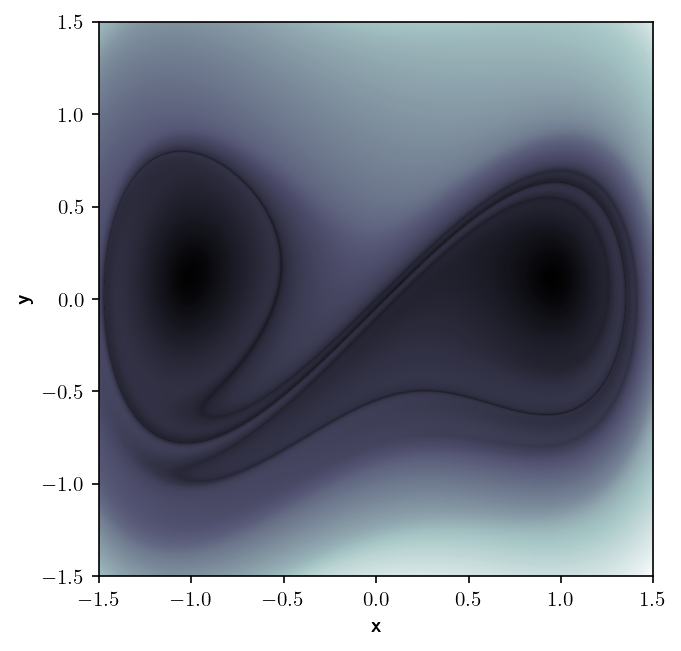} b)\includegraphics[scale=0.5]{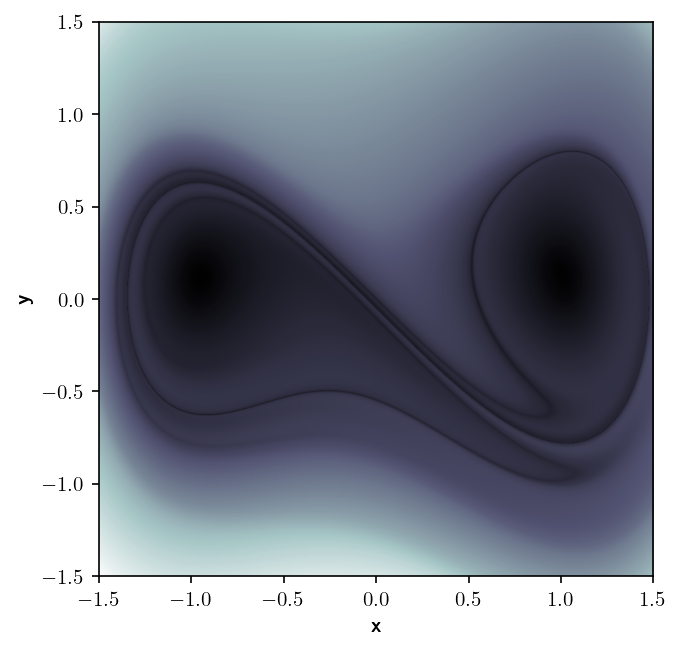} 
  \end{center}
  \caption{ a) Evaluation of the backwards integration of \eqref{metricMd} on the Duffing equation \eqref{duff} at $t_0=0$ and $\tau=10$. This represents the unstable manifold; b) Evaluation of the forwards integration of \eqref{metricMd} on the Duffing equation \eqref{duff} at $t_0=0$ and $\tau=10$. This represents the stable manifold.    }
  \label{fig:duffspl}
\end{figure}

This section discusses further results on the Lagrangian Uncertainty Quantifier by 
considering the evaluation of
 \eqref{eq:UQ1} over a vector field obtained from the nonlinear periodically forced Duffing equation: 
\begin{eqnarray}
 \dot{x} &=& y \nonumber \\
  \dot{y} &=& x-x^3+ \epsilon \sin t  
  \label{duff}
\end{eqnarray}
Prior to discuss outputs of \eqref{eq:UQ1} into Eq. \eqref{duff}, we discuss the structure of invariant manifolds 
of hyperbolic trajectories in Eq. \eqref{duff} for the case $\epsilon=0.1$ and the persistence versus this time dependent perturbation of tori present in the unforced version of Eq. \eqref{duff}, i. e, $\epsilon=0$. In the perturbed case, the hyperbolic fixed point placed at the origin, becomes a hyperbolic periodic trajectory \cite{madrid2009}, and their stable and unstable manifold can be highlighted by the Equation \ref{metricMd} appeared in Section 2 \cite{mancho2013}. Additionally, as discussed in \cite{lopesino2017theoretical} a scale factor $1/(2\tau)$ applied to \eqref{metricMd} converts the expression to an average, which in compact Hamiltonian systems like \eqref{duff}  converges for $\tau \to \infty$ and when this convergence is observed, level curves correspond to invariant structures of the dynamical system. Convergence of means are computationally verifiable on tori, however on hyperbolic sets as discussed in  \cite{garcia2018detection}, rounding computational errors practically prevent convergence.
\begin{figure}[htb!]
  \begin{center}
 a) \includegraphics[scale=0.5]{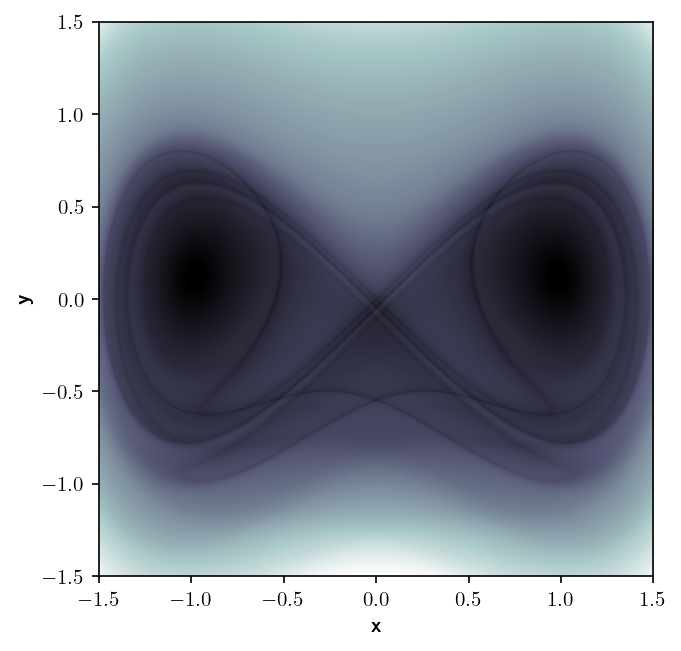} b)\includegraphics[scale=0.5]{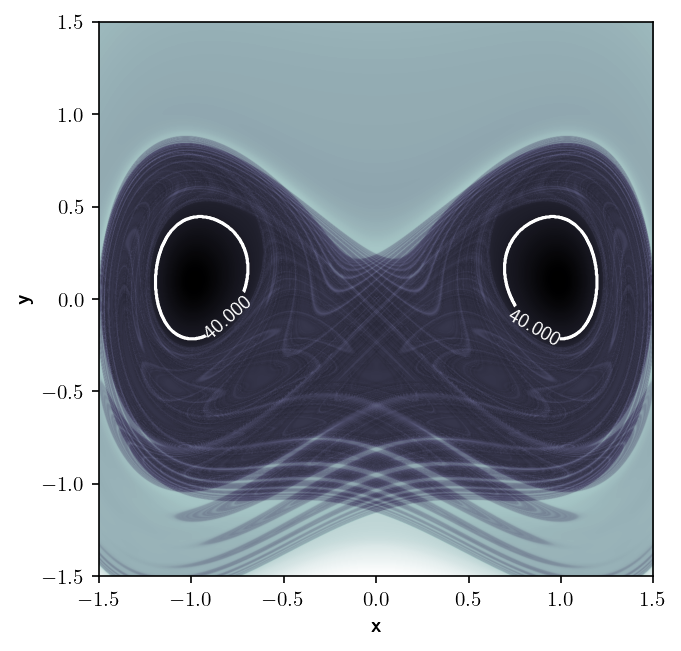} 
  \end{center}
  \caption{ Evaluation of \eqref{metricMd} for the Duffing equation \eqref{duff} at $t_0=0$. a) $\tau=10$; b) $\tau=50$. At this $\tau$ the average of $M$ in the smooth regions has converged and the level curve represent and invariant tori.    }
  \label{fig:duff}
\end{figure}
Figure \ref{fig:duff} illustrates these points.  Panel a) shows invariant manifolds related to the hyperbolic trajectory that are obtained from  \eqref{metricMd} for $\tau=10$. These structures become much more rich for larger $\tau$ as panel b) confirms. In this panel since the average of $M$ has converged in the smooth region, level curves in the area highlight tori.
The outputs of Figure \ref{fig:duff}a) may be split into two figures, by separating  the backwards and forwards integration of Eq.  \eqref{metricMd}, which displays, respectively, the unstable and stable manifolds. These outputs are depicted, respectively, in Figure \ref{fig:duffspl}a) and b).

\begin{figure}[htb!]
 \begin{center}
  a)\includegraphics[scale=0.5]{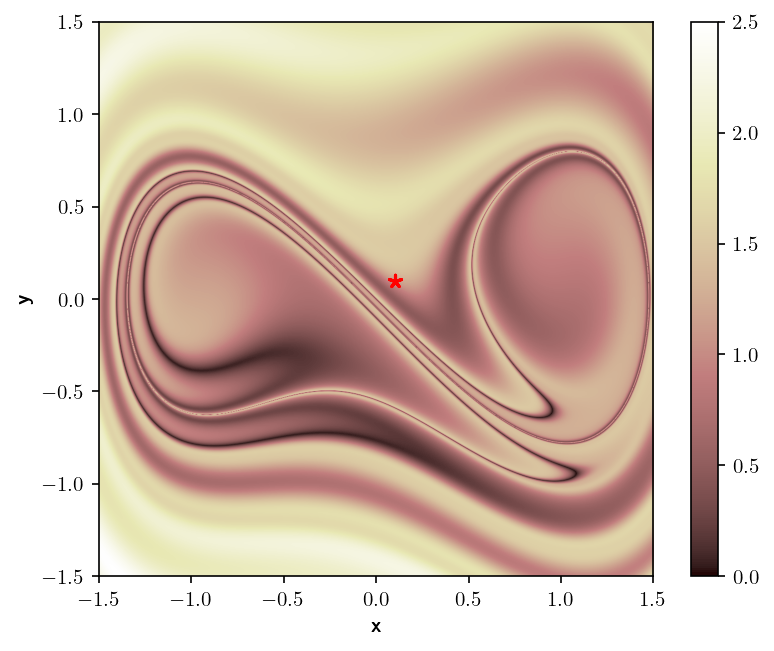}
  b) \includegraphics[scale=0.5]{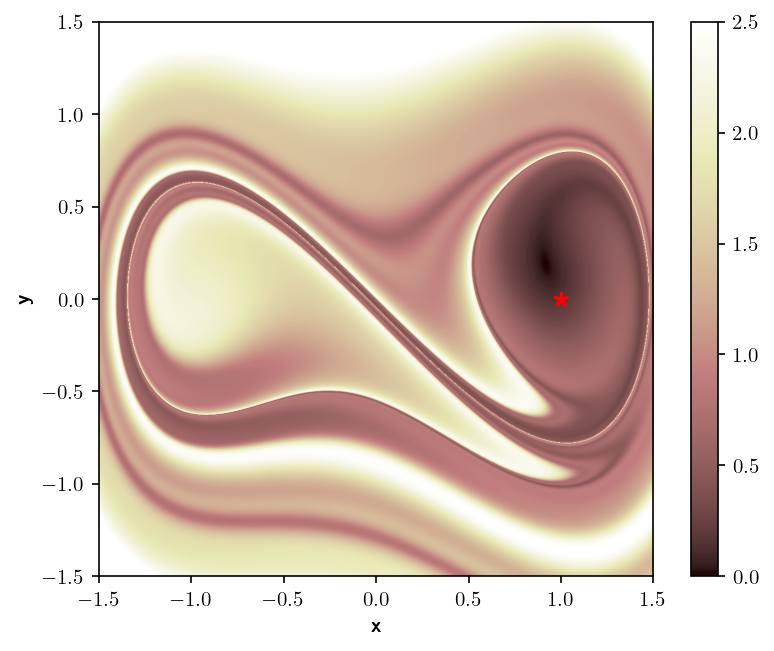} 
 c)\includegraphics[scale=0.5]{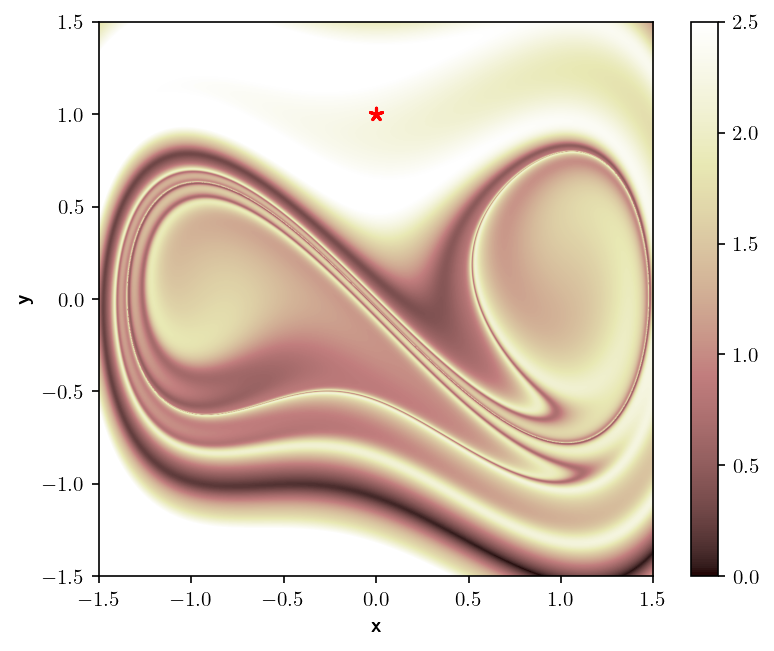} 
\end{center}
   \caption{ Evaluation of Eq.\eqref{eq:UQ1} 
 on the Duffing equation \eqref{duff} at $t_0=0$, $t=t_0+\tau=10$ with different targets. a) Target ${\bf x}^*=(0.1, 0.1)$; b) target ${\bf x}^*=(1, 0)$; c) target ${\bf x}^*=(0, 1)$. }
   \label{fig:duffuqld}
 \end{figure}

Figure \ref{fig:duffuqld} shows the evaluation of $L_{UQ}$ as in  Eq.\eqref{eq:UQ1} with $p=2$ for $t_0=0$, $t=t_0+\tau=10$ and different targets. In panel a) the chosen target is ${\bf x}^*=(0.1, 0.1)$. This target is
within the chaotic region displayed in Fig. \ref{fig:duff}b), close to the unstable manifold of the Duffing equation. It is observed that minimum values are reached on the stable manifold and that the structure of $L_{UQ}$ is correlated to it.  This setting is similar to what was observed for the ocean case described in Section 2, in which the observed oil evolution is aligned with the unstable manifold and minimum values of $L_{UQ}$ are found on the stable manifold that are optimal pathways towards the unstable manifolds.  In panel b) we show the results for the
target  ${\bf x}^*=(1, 0)$ which is within the right tori like structure  highlighted  in Fig. \ref{fig:duff}b). Accordingly, uncertainty values are very low for initial observations $(x_0, y_0)$ in the corresponding tori region, but very high for the tori like region at the left side. Indeed initial observations in this region will never go near a  final observation in the right tori like structure, and therefore this model is {\em structurally uncertain} for those observations, i.e. the model is inadequate to represent those.  Panel c) shows the computation of the uncertainty for a target ${\bf x}^*=(0,1)$ outside the geometry of the unstable invariant manifold displayed in Fig. \ref{fig:duffspl}a). It is remarkable the persistence in all these examples of an structure on the uncertainty field with singular features linked to the stable manifold  independently of the target value    ${\bf x}^*$.
\begin{figure}[htb!]
\begin{center}
 a) \includegraphics[scale=0.3]{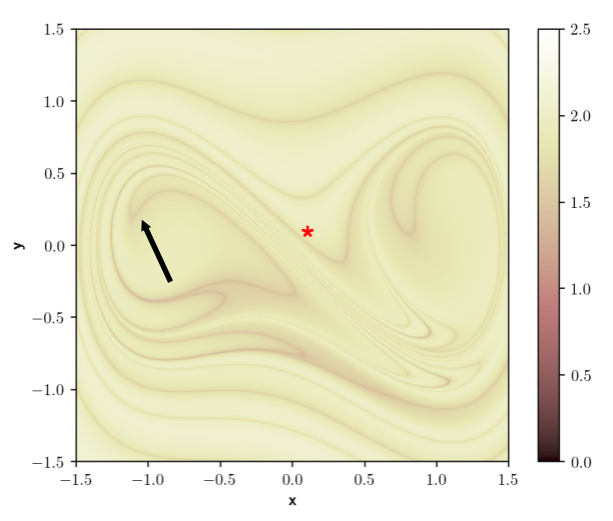} 
 b)\includegraphics[scale=0.32]{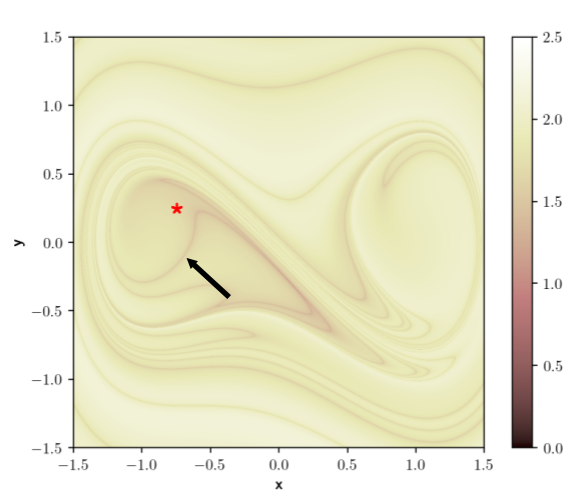} 
 c)\includegraphics[scale=0.3]{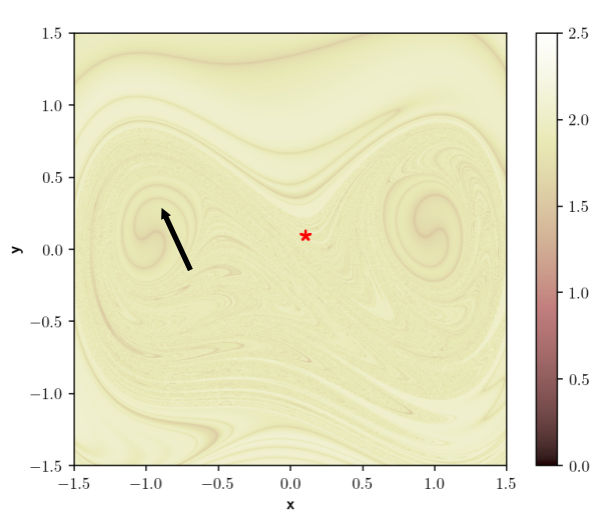} 
 \end{center}
   \caption{ Evaluation of Eq.\eqref{eq:UQ2} 
 on the Duffing equation \eqref{duff} at $t_0=0$, for $p= 0.1$
 and  different targets. Black arrows highlight features that for this system are known to be spurious. a) $t_0=0$, $t=t_0+\tau=10$ and target ${\bf x}^*=(0.1, 0.1)$; b) $t_0=0$, $t=t_0+\tau=10$ and target ${\bf x}^*=(-0.75, 0.25)$; c)  $t_0=0$, $t=t_0+\tau=50$ and target ${\bf x}^*=(0.1, 0.1)$  }
  \label{fig:duffuqsp}
\end{figure}


Figure \ref{fig:duffuqsp} shows the evaluation of  Eq.\eqref{eq:UQ2} with $p=0.1$, and different targets or integration periods. Panels a) and b) have the same integration periods, but different targets.  In them the black arrow points out different spurious features that seem singular but {\it do not correspond to any invariant manifold}. Figure \ref{fig:duffuqsp} c) shows the same than panel a) for an integration period of 50. The black arrow marks a region, which from the analysis of Fig. \ref{fig:duff}, is known to be covered by tori, while  the structure attained from Eq.\eqref{eq:UQ2}  does not highlight this.
Indeed, the ergodic partition theory discussed in \cite{mezic99,susuki2004} and implemented in \cite{lopesino2017theoretical} for Lagrangian Descriptors such as that in Eq. \eqref{metricMd}, shows that this ability to highlight tori, requires averaging along trajectories, while expressions  Eq.  \eqref{eq:UQ1} or  \eqref{eq:UQ2} {\it are not averages}, and therefore do not qualify for the application of ergodic principles. These expressions are not suitable for identifying this kind of invariant objects.


\section{Discussion}

There is much interest in uncertainty quantification involving trajectories in ocean data sets.
Recent efforts in this direction are \cite{kirwan2003,huntley2011,feng2019, rypina20}.  Vieira et al. \cite{rypina20} have developed a clustering method to partition the space of trajectory data sets into distinct flow regions. This method contains free parameters and they use uncertainty quantification to assess the parameter dependence on the partitions they obtain. The work \cite{ feng2019} discusses  approaches for uncertainty
quantification from a geometrical perspective, and  uncertainty quantifiers based on distances are proposed, but no connections are proposed with invariant dynamical structures. Similarly in \cite{huntley2011} uncertainty quantifiers based on a descriptive statistics of distances between modelled trajectories and observations  are used, but no relations are presented between these and the invariant dynamical structures.
Results  appeared in \cite{haller2020} also discuss along these lines, but in the context of general models, not related to  ocean data. In this work authors have established  
links between uncertainty quantification and  invariant manifolds  under explicit mathematical assumptions for the "true" model. The results in our work  do not require  such assumptions. 

Equations similar to \eqref{eq:UQ1} or  \eqref{eq:UQ2} have been used in the literature to highlight Lagrangian structures in oceanic flows.
For instance  \cite{prants2011numerical,prants2014lagrangian}  have done so.  Prants in  \cite{prants2011numerical} proposed to use the arc-length $D$ to study the displacement of particles in the coasts of Japan. 

\begin{equation}
D = \sqrt{(x_f-x_0)^2+(y_f-y_0)^2}. \label{prants}
\end{equation}

\noindent
Here, $D$ represents the relative displacement  of  a  particle  from 
its  initial  position  $(x_0, y_0)$  to certain  final  one  $(x_f, y_f)$. This expression is the analogue to Eq.\eqref{eq:UQ1} with $p=2$ and target adjusted to each initial condition.  However, in this work no connections are established between Eq. \eqref{prants} and  uncertainty quantification.


\section{Conclusions}

This article explores the implications of a definition for Uncertainty Quantification recently proposed in oceanic contexts \cite{garciasanchez2020}. It is found that for this definition, which is associated to forward  Uncertainty Quantification, stable invariant manifolds of hyperbolic trajectories of the underlying flow, provide a structure for it. That is, we found that the proposed Uncertainty Quantifier is a function that contains  a very rich structure, which  is related to these well known structures from dynamical systems theory. For selected examples this statement has been proven. Also examples are discussed in which singular structures of the Uncertainty Quantifier highlight spurious structures,  not aligned with invariant features of the dynamical system. Other invariant structures like tori, are not recovered by the Uncertainty Quantifier. Also, connections between UQ and invariant dynamical structures have been used to provide a framework for discussion on structural uncertainty, which is related to inadequate models.  
 This vision enriches traditional descriptions of UQ on which structure is discussed just in terms of means and/or statistical moments of distributions.

The findings described in  this article are particularly interesting because they have important environmental applications. Nowadays, multiple ocean data sources are available   and in this context our results allow  a quantitative comparison  of the transport properties associated to them. In this way, discriminating the level of performance of different data sources will help to gain precision in the description of dispersion of contaminants, determination of waste and plastic sources, etc.

\section*{Acknowledgements}
GGS and AMM acknowledge support from IMPRESSIVE, a project funded by the European Union's Horizon 2020 research and innovation programme under grant agreement No 821922. SW acknowledges the support of ONR Grant No.~N00014-01-1-0769.

\bibliographystyle{ieeetr}
\bibliography{sciref,LD}{}

\begin{thebibliography}{10}

\bibitem{sullivan2015}
T.~J. Sullivan, {\em Introduction to Uncertainty Quantification}.
\newblock Springer, 2015.

\bibitem{garciasanchez2020}
G.~Garc{\'\i}a-S{\'a}nchez, A.~M. Mancho, A.~G. Ramos, J.~Coca,
  B.~P{\'e}rez-G{\'o}mez, E.~{\'A}lvarez-Fanjul, M.~G. Sotillo,
  M.~Garc{\'\i}a-Le{\'o}n, V.~J. Garc{\'\i}a-Garrido, and S.~Wiggins, ``Very
  high resolution tools for the monitoring and assessment of environmental
  hazards in coastal areas,'' {\em Frontiers in Marine Science}, vol.~7,
  no.~605804, 2021.

\bibitem{sotillo2020coastal}
M.~G. Sotillo, P.~Cerralbo, P.~Lorente, M.~Grifoll, M.~Espino,
  A.~Sanchez-Arcilla, and E.~{\'A}lvarez-Fanjul, ``Coastal ocean forecasting in
  {S}panish ports: the {SAMOA} operational service,'' {\em Journal of
  Operational Oceanography}, vol.~13, no.~1, pp.~37--54, 2020.

\bibitem{madrid2009}
J.~A.~J. Madrid and A.~M. Mancho, ``Distinguished trajectories in time
  dependent vector fields,'' {\em Chaos}, vol.~19, p.~013111, 2009.

\bibitem{mendoza2010}
C.~Mendoza and A.~M. Mancho, ``The hidden geometry of ocean flows,'' {\em Phys.
  Rev. Lett.}, vol.~105, no.~3, p.~038501, 2010.

\bibitem{mancho2013}
A.~M. Mancho, S.~Wiggins, J.~Curbelo, and C.~Mendoza, ``Lagrangian descriptors:
  A method for revealing phase space structures of general time dependent
  dynamical systems,'' {\em Communications in Nonlinear Science and Numerical
  Simulations}, vol.~18, no.~12, pp.~3530--3557, 2013.

\bibitem{lopesino2017theoretical}
C.~Lopesino, F.~Balibrea-Iniesta, V.~J. Garc{\'\i}a-Garrido, S.~Wiggins, and
  A.~M. Mancho, ``A theoretical framework for lagrangian descriptors,'' {\em
  International Journal of Bifurcation and Chaos}, vol.~27, no.~01, p.~1730001,
  2017.

\bibitem{lopesino2015lagrangian}
C.~Lopesino, F.~Balibrea, S.~Wiggins, and A.~M. Mancho, ``Lagrangian
  descriptors for two dimensional, area preserving, autonomous and
  nonautonomous maps,'' {\em Communications in Nonlinear Science and Numerical
  Simulation}, vol.~27, no.~1-3, pp.~40--51, 2015.

\bibitem{garcia2018detection}
V.~J. Garc{\'\i}a-Garrido, F.~Balibrea-Iniesta, S.~Wiggins, A.~M. Mancho, and
  C.~Lopesino, ``Detection of phase space structures of the cat map with
  lagrangian descriptors,'' {\em Regular and Chaotic Dynamics}, vol.~23, no.~6,
  pp.~751--766, 2018.

\bibitem{moser56}
J.~Moser, ``The analytic invariants of an area-preserving mapping near a
  hyperbolic fixed point,'' {\em Comm. Pure App. Math.}, vol.~9, pp.~673--692,
  1956.

\bibitem{mezic99}
I.~Mezic and S.~Wiggins, ``A method for visualization of invariant sets of
  dynamical systems based on the ergodic partition,'' {\em Chaos}, vol.~9,
  pp.~213--218, 1999.

\bibitem{susuki2004}
Y.~Susuki and I.~Mezic, ``Ergodic partition of phase space in continuous
  dynamical systems,'' {\em Joint 48th IEEE Conf. Decision and Control and 28th
  Chinese Control Conf.}, pp.~7497--7502, 204.

\bibitem{kirwan2003}
A.~D. Kirwan~Jr., M.~Toner, and L.~Kantha, ``Predictability, uncertainty, and
  hyperbolicity in the ocean,'' {\em International Journal of Engineering
  Science}, vol.~41, p.~249, 2003.

\bibitem{huntley2011}
S.~H. Huntley, B.~L. Lipphardt~Jr., and A.~D. Kirwan~Jr., ``Lagrangian
  predictability assessed in the east china sea,'' {\em Ocean Modelling},
  vol.~36, p.~163, 2011.

\bibitem{feng2019}
D.~Feng, P.~Passalacqua, and B.~R. Hodges, ``Innovative approaches for
  geometric uncertainty quantification in an operational oil spill modeling
  system,'' {\em Journal of Marine Science and Engineering}, vol.~7, no.~8,
  p.~259, 2019.

\bibitem{rypina20}
G.~S. Vieira, I.~I. Rypina, and M.~R. Allshouse, ``Uncertainty quantification
  of trajectory clustering applied to ocean ensemble forecasts,'' {\em Fluids},
  vol.~5, no.~4, p.~184, 2020.

\bibitem{haller2020}
B.~Kasz\'as and G.~Haller, ``Universal upper estimate for prediction errors
  under moderate model uncertainty,'' {\em Chaos}, vol.~30, p.~113144, 2020.

\bibitem{prants2011numerical}
S.~Prants, M.~Y. Uleysky, and M.~Budyansky, ``Numerical simulation of
  propagation of radioactive pollution in the ocean from the fukushima dai-ichi
  nuclear power plant,'' in {\em Doklady Earth Sciences}, vol.~439, p.~1179,
  Springer, 2011.

\bibitem{prants2014lagrangian}
S.~Prants, M.~Budyansky, and M.~Y. Uleysky, ``Lagrangian study of surface
  transport in the kuroshio extension area based on simulation of propagation
  of fukushima-derived radionuclides,'' {\em Nonlinear Processes in
  Geophysics}, vol.~21, no.~1, pp.~279--289, 2014.

\end{thebibliography}

\end{document}